\pdfoutput=1

\documentclass[11pt]{amsart}
\usepackage{amssymb,latexsym,amsmath}
\usepackage[pdftex]{graphicx}
\usepackage{graphics}
\usepackage{harvard,enumerate,amssymb}
\usepackage{blkarray}
\usepackage{bm}

\newcommand{\indep}{{\bot\negthickspace\negthickspace\bot}}

\newtheorem{thm}{Theorem}[section]

\newtheorem{lem}{Lemma}[section]

\newtheorem{ass}{Assumption}[section]

\newtheorem{condition}{Condition}[section]

\theoremstyle{definition}
\newtheorem{defn}{Definition}[section]
\newtheorem{example}{Example}[section]
\newtheorem{rem}{Remark}[section]
\newtheorem*{notation}{Notation}
\numberwithin{equation}{section}

\newcommand{\argmin}{\operatornamewithlimits{argmin}}

\begin{document}
\title[]{Nonparametric Analysis of Random Utility Models}
\author[Kitamura]{Yuichi Kitamura$^*$}
\address{Cowles Foundation for Research in Economics, Yale University, New
Haven, CT 06520.}
\email{yuichi.kitamura@yale.edu}
\author[Stoye]{J\"org Stoye$^{**}$}
\address{Departments of Economics, Cornell University, Ithaca, NY 14853, and University of Bonn, Germany; Hausdorff Center for Mathematics, Bonn, Germany.}
\email{stoye@cornell.edu}
\date{This Version: May 16, 2018.}
\thanks{\emph{Keywords}: Stochastic Rationality}
\thanks{JEL Classification Number: C14}
\thanks{$^{\ast }$ $^{\ast \ast }$ Kitamura acknowledges financial support
from the National Science Foundation via grants SES-0851759 and SES-1156266.
Stoye acknowledges support from the National Science Foundation under grant
SES-1260980 and through Cornell University's Economics Computer Cluster Organization, which was partially funded through NSF Grant SES-0922005. We thank Donald Andrews, Richard Blundell, Rahul Deb, Phil Haile, Chuck Manski, Daniel Martin, Rosa Matzkin, Francesca Molinari, participants at various seminars, and in particular, Elie Tamer and numerous referees for helpful comments.  We are grateful to Richard Blundell and Ian Crawford for their help with the FES data.  We would like to thank Whitney Newey in particular for suggesting the use of a control function to account for endogeneity in total expenditure.  All errors are our own.  Ben Friedrich, Duk Gyoo Kim, Matthew Thirkettle, and Yang Zhang provided excellent research assistance.}

\begin{abstract}
This paper develops and implements a nonparametric test of Random Utility Models. The motivating application is to test the null hypothesis that a sample of cross-sectional demand distributions was generated by a population of rational consumers. We test a necessary and sufficient condition for this that does not restrict unobserved heterogeneity or the number of goods. We also propose and implement a control function approach to account for endogenous expenditure. An econometric result of independent interest is a test for linear inequality constraints when these are represented as the vertices of a polyhedral cone rather than its faces. An empirical application to the U.K. Household Expenditure Survey illustrates computational feasibility of the method in demand problems with 5 goods.
\end{abstract}

\maketitle

\section{Introduction}

\label{sec:intro}This paper develops new tools for the nonparametric analysis of Random Utility Models (RUM). We test the null hypothesis that a repeated cross-section of demand data might have been generated by a population of rational consumers, without restricting either unobserved heterogeneity or the number of goods. Equivalently, we empirically test McFadden and Richter's (1991) \textit{Axiom of Revealed Stochastic Preference}. To do so, we develop a new statistical test that promises to be useful well beyond the motivating application.

We start from first principles and end with an empirical application. Core contributions made along the way are as follows.

First, the testing problem appears formidable: A structural parameterization of the null hypothesis would involve an essentially unrestricted distribution over all nonsatiated utility functions. However, the problem can, without loss of information, be rewritten as one in which the universal choice set is finite. Intuitively, this is because a RUM only restricts the population proportions with which preferences between different budgets are directly revealed. The corresponding sample information can be preserved in an appropriate discretization of consumption space.

More specifically, observable choice proportions must be in the convex hull of a finite (but long) list of vectors. Intuitively, these vectors characterize rationalizable nonstochastic choice types, and observable choice proportions are a mixture over them that corresponds to the population distribution of types. This builds on \citeasnoun{mcfadden-2005} but with an innovation that is crucial for testing: While the set just described is a finite polytope, the null hypothesis can be written as a cone. Furthermore, computing the list of vectors is hard, but we provide algorithms to do so efficiently.

Next, the statistical problem is to test whether an estimated vector of choice proportions is inside a nonstochastic, finite polyhedral cone. This is reminiscient of multiple linear inequality testing and shares with it the difficulty that inference must take account of many nuisance parameters. However, in our setting, inequalities are characterized only implicitly through the vertices of their intersection cone. It is not computationally possible to make this characterization explicit. We provide a novel test and prove that it controls size uniformly over a reasonable class of d.g.p.'s without either computing facets of the cone or resorting to globally conservative approximation. This is a contribution of independent interest that has already seen other applications \cite{DKQS16,Hubner,LQS15,LQS18}. Also, while our approach can become computationally costly in high dimensions, it avoids a statistical curse of dimensionality (i.e., rates of approximation do not deteriorate), and our empirical exercise shows that it is practically applicable to at least five-dimensional commodity spaces.    

Finally, we leverage recent results on control functions (\citeasnoun{imbens2009}; see also \citeasnoun{blundell2003}) to deal with endogeneity for unobserved heterogeneity of unrestricted dimension. These contributions are illustrated on the U.K. Family Expenditure Survey, one of the work horse data sets of the literature. In that data, estimated demand distributions are not stochastically rationalizable, but the rejection is not statistically significant.

The remainder of this paper is organized as follows. Section \ref{sec:literature} discusses the related literature. Section \ref{sec:population} lays out the model, develops a geometric characterization of its empirical content, and presents algorithms that allow one to compute this characterization in practice. All of this happens at population level, i.e. all identifiable quantities are known. Section \ref{sec:Inference} explains our test and its implementation under the assumption that one has an estimator of demand distributions and an approximation of its sampling distribution. Section \ref{sec:extending} explains how to get the estimator, and a bootstrap approximation to its distribution, by both smoothing over expenditure and adjusting for endogenous expenditure. Section \ref{sec:monte carlo} contains a Monte Carlo investigation of the test's finite sample performance, and Section \ref{sec:empirics} contains our empirical application. Section \ref{sec:conclusion} concludes. Supplemental materials collect all proofs (Appendix A), pseudocode for some algorithms (Appendix B), and some algebraic elaborations (Appendix C).  

\section{Related Literature}

\label{sec:literature}

Our framework for testing Random Utility Models is built from scratch in the sense that it only presupposes classic results on nonstochastic revealed preference, notably the characterization of individual level rationalizability through the Weak \cite{Samuelson38}, Strong \cite{Houthakker50}, or Generalized \cite{Afriat67} Axiom of Revealed Preference (WARP, SARP, and GARP henceforth). At the population level, stochastic rationalizability was analyzed in classic work by \citeasnoun{mcfadden-richter} updated by \citeasnoun{mcfadden-2005}. This work was an important inspiration for ours, and we will further clarify the relation later, but they did not consider statistical testing nor attempt to make the test operational, and could not have done so with computational constraints even of 2005.

An influential related research project is embodied in a sequence of papers by Blundell, Browning, and Crawford (2003, 2007, 2008; BBC henceforth), where the 2003 paper focuses on testing rationality and bounding welfare and later papers focus on bounding counterfactual demand. BBC assume the same observables as we do and apply their method to the same data, but they analyze a nonstochastic demand system generated by nonparametric estimation of Engel curves. This could be loosely characterized as revealed preference analysis of a representative consumer and in practice of average demand. \citeasnoun{Lewbel01} gives conditions on a RUM that ensure integrability of average demand, so BBC effectively add those assumptions to ours. Also, the nonparametric estimation step in practice constrains the dimension of commodity space, which equals three in their empirical applications.\footnote{BBC's implementation exploits only WARP and therefore a necessary but not sufficient condition for rationalizability. This is remedied in \citeasnoun{BBCCDV}.} 

\citeasnoun{Manski07} analyzes stochastic choice from subsets of an abstract, finite choice universe. He states the testing and extrapolation problems in the abstract, solves them explicitly in simple examples, and outlines an approach to non-asymptotic inference. (He also considers models with more structure.) While we start from a continuous problem and build a (uniform) asymptotic theory, the settings become similar after our initial discretization step. However, methods in \citeasnoun{Manski07} will only be practical for choice universes with a handful of elements, an order of magnitude less than in Section \ref{sec:empirics} below. In a related paper, \citeasnoun{Manski14} uses our computational toolkit for choice extrapolation.

Our setting much simplifies if there are only two goods, an interesting but obviously very specific case. \citeasnoun{BKM14} bound counterfactual demand in this setting through bounding quantile demands. They justify this through an invertibility assumption. \citeasnoun{HS15} show that with two goods, this assumption has no observational implications.\footnote{A similar point is made, and exploited, by \citeasnoun{HN16}.} Hence, \citeasnoun{BKM14} use the same assumptions as we do; however, the restriction to two goods is fundamental. \citeasnoun{BKM17} conceptually extend this approach to many goods, in which case invertibility is a restriction. A nonparametric estimation step again limits the dimensionality of commodity space. They apply the method to similar data and the same goods as BBC, meaning that their nonparametric estimation problem is two-dimensional.

\citeasnoun{HN16} nonparametrically bound average welfare under assumptions resembling ours, though their approach additionally imposes smoothness restrictions to facilitate nonparametric estimation and interpolation.  Their main identification results apply to an arbitrary number of goods, but the approach is based on nonparametric smoothing, hence the curse of dimensionality needs to be addressed. The empirical application is to two goods. 

With more than two goods, pairwise testing of a stochastic analog of WARP amounts to testing a necessary but not sufficient condition for stochastic rationalizability. This is explored by \citeasnoun{HS14} in a setting that is otherwise ours and also on the same data. \citeasnoun{Kawaguchi} tests a logically intermediate condition, again on the same data. A different test of necessary conditions was proposed by \citeasnoun{Hoderlein11}, who shows that certain features of rationalizable individual demand, like adding up and standard properties of the Slutsky matrix, are inherited by average demand under weak conditions. The resulting test is passed by the same data that we use. \citeasnoun{DHN16} propose a similar test using quantiles.

Section \ref{sec:Inference} of this paper is (implicitly) about testing multiple inequalities, the subject of a large literature in economics and statistics. See, in particular, \citeasnoun{ghm} and \citeasnoun{wolak-1991} and also \citeasnoun{chernoff}, \citeasnoun{kudo63}, \citeasnoun{perlman}, \citeasnoun{shapiro}, and \citeasnoun{takemura-kuriki} as well as \citeasnoun{andrews-hac}, \citeasnoun{BCS15}, and \citeasnoun{guggenberger-hahn-kim}. For the also related setting of inference on parameters defined by moment inequalities, see furthermore \citeasnoun{andrews-soares}, \citeasnoun{bugni-2010}, \citeasnoun{Canay10}, \citeasnoun{cht}, \citeasnoun{imbens-manski}, \citeasnoun{romano-shaikh}, \citeasnoun{rosen-2008}, and \citeasnoun{stoye-2009}.  The major difference to these literatures is that moment inequalities, if linear (which most of the papers do not assume), define a polyhedron through its faces, while the restrictions generated by our model correspond to its vertices. One cannot in practice switch between these representations in high dimensions, so that we have to develop a new approach. This problem also occurs in a related problem in psychology, namely testing if binary choice probabilities are in the so-called Linear Order Polytope. Here, the problem of computing explicit moment inequalities is researched but unresolved (e.g., see \citeasnoun{DR15} and references therein), and  we believe that our test is of interest for that literature. Finally, \citeasnoun{HS14} only compare two budgets at a time, and \citeasnoun{Kawaguchi} tests necessary conditions that are directly expressed as moment inequalities. Therefore, inference in both papers is much closer to the aforecited literature.

\section{Analysis of Population Level Problem}

\label{sec:population}

We now show how to verify rationalizability of a known set of cross-sectional demand distributions on $J$ budgets. The main results are a tractable geometric characterization of stochastic rationalizability and algorithms for its practical implementation.

\subsection{Setting up the Model.}
\label{sec:setup}

Throughout this paper, we assume the existence of $J < \infty$ fixed budgets $\mathcal{B}_j$ characterized by price vectors $p_j \in \mathbf{R}^K_+$ and expenditure levels $W_j >0$. Normalizing $W_j=1$ for now, we can write these budgets as
\begin{equation*}
\mathcal{B}_{j}=\{y \in \mathbf{R}^K_+:p_j'y=1\},j=1,...,J.
\end{equation*}%
We also start by assuming that the corresponding cross-sectional distributions of demand are known. Thus, assume that demand in budget $\mathcal{B}_j$ is described by the random variable $y(p_j)$, then we know  
\begin{equation}
\label{eq:define_pi_j}
P_j(x):=\Pr(y(p_j)\in x), \ x\subset \mathbf{R}_{+}^{K}
\end{equation}%
for $j=1,...,J$.\footnote{To keep the presentation simple, here and henceforth we are informal about probability spaces and measurability. See \citeasnoun{mcfadden-2005} for a formally rigorous setup.} We will henceforth call $(P_1,...,P_J)$ a \textit{stochastic demand system}.

The question is if this system is rationalizable by a RUM. To define the latter, let
\begin{equation*}
u:\mathcal{\mathbf{R}}_+^K \mapsto \mathbf{R}
\end{equation*}%
denote a utility function over consumption vectors $y \in \mathbf{R}_+^K$. Consider for a moment an individual consumer endowed with some fixed $u$, then her choice from a budget characterized by normalized price vector $p$ would be
\begin{equation} \label{eq:intro}
y \in \arg \max_{y\in \mathbf{R}_{+}^{K}:p'y \leq 1}u(y),
\end{equation}
with arbitrary tie-breaking if the solution is not unique. For simplicity, we restrict utility functions by monotonicity (\textquotedblleft more is better\textquotedblright) so that choice is on budget planes, but this is not conceptually necessary.

The RUM postulates that
$$
u \sim P_u,
$$
i.e. $u$ is not constant but is distributed according to a constant (in $j$) probability law $P_u$. In our motivating application, $P_u$ describes the distribution of preferences in a population of consumers, but other interpretations are conceivable. For each $j$, the random variable $y(p_j)$ then is the distribution of $y$ defined in \eqref{eq:intro} that is induced by $p_j$ and $P_u$. Formally:

\begin{defn}
\label{def:model}
The stochastic demand system $(P_1,...,P_J)$ is \textit{(stochastically) rationalizable} if there exists a distribution $P_u$ over utility functions $u$ so that
\begin{equation}
P_j (x)=\int \! 1\{\arg\max_{y\in \mathbf{R}_+^K:p_j'y=1}u(y) \in x\} \mathrm{d}P_u, \quad x\subset \mathcal{B}_j,j=1,...,J. \label{eq:model}
\end{equation}
\end{defn}

This model is completely parameterized by $P_u$, but it only partially identifies $P_u$ because many distinct $P_u$ will induce the same stochastic demand system. We do not place substantive restrictions on $P_u$, thus we allow for minimally constrained, infinite dimensional unobserved heterogeneity across consumers.

Definition \ref{def:model} reflects some simplifications that we will drop later. First, $W_j$ and $p_j$ are nonrandom, which is the framework of \citeasnoun{mcfadden-richter} and others but may not be realistic in applications. In the econometric analysis in Section \ref{sec:extending} as well as in our empirical analysis in Section \ref{sec:Empirical Application}, we treat $W_j$ as a random variable that may furthermore covary with $u$. Also, we initially assume that $P_u$ is the same across price regimes. Once $W_j$ (hence $p_j$, after income normalization) is a random variable, this is essentially the same as imposing $W_j \indep u$, an assumption we maintain in Section \ref{sec:Inference} but drop in Section \ref{sec:extending} and in our empirical application. However, for all of these extensions, our strategy will be to effectively reduce them to \eqref{eq:model}, so testing this model is at the heart of our contribution.

\subsection{A Geometric Characterization.}
\label{sec:geometry}

The model embodied in \eqref{eq:model} is extremely general; again, a parameterization would involve an essentially unrestricted distribution over utility functions. However, we next develop a simple geometric characterization of the model's empirical content and hence of stochastic rationalizability.

To get an intuition, consider the simplest example in which \eqref{eq:model} can be tested:

\begin{example}
\label{example 1}
There are two intersecting budgets, thus $J=2$ and there exists $y \in \mathbf{R}^K_{++}$ with $p_1'y=p_2'y$. 
\end{example}

\begin{figure}
\centering
\includegraphics[scale=1]{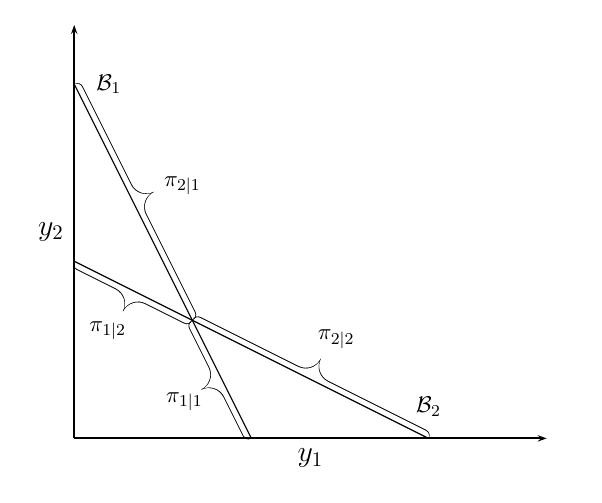}
\caption{Visualization of Example \ref{example 1}.}
\label{figure 1}
\end{figure}    

Consider Figure \ref{figure 1}, whose labels will become clear. (The restriction to $\mathbf{R}^2$ is only for the figure.) It is well known that in this example, individual choice behavior is rationalizable unless choice from each budget is below the other budget, in which case a consumer would revealed prefer each budget to the other one. Does this restrict repeated cross-section choice probabilities? Yes: Supposing for simplicity that there is no probability mass on the intersection of budget planes, it is easy to see (e.g. by applying Fr\'{e}chet-Hoeffding bounds) that the cross-sectional probabilities of the two line segments labeled $(\pi_{1|1},\pi_{1|2})$ must not sum to more than $1$. This condition is also sufficient \cite{Matzkin06}.

Things rapidly get complicated as budgets are added, but the basic insight scales. The only relevant information for testing \eqref{eq:model} is what fractions of consumers revealed prefer budget $j$ to $k$ for different $(j,k)$. This information is contained in the cross-sectional choice probabilities of the line segments highlighted in Figure \ref{figure 1} (plus, for noncontinuous demand, the intersection). The picture will be much more involved in interesting applications -- see Figure \ref{figure 2} for an example -- but the idea remains the same. This insight allows one to replace the universal choice set $\mathbf{R}^K_+$ with a finite set and stochastic demand systems with lists of corresponding choice probabilities. But then there are only finitely many rationalizable nonstochastic cross-budget choice patterns. A rationalizable stochastic demand system must be a mixture over them and therefore lie inside a certain finite polytope.

Formalizing this insight requires some notation.

\begin{defn}\label{def:patches}
Let $\mathcal{X} :=\{x_1,...,x_I\}$ be the coarsest partition of $\cup _{j=1}^J\mathcal{B}_j$ such that for any $i \in \{1,...,I\}$ and $j \in \{1,...,J\}$, $x_i$ is either completely on, completely strictly above, or completely strictly below budget plane $\mathcal{B}_j$. Equivalently, any $y_1,y_2 \in \cup _{j=1}^J\mathcal{B}_j$ are in the same element of the partition iff $\text{sg}(p_j'y_1-1)=\text{sg}(p_j'y_2-1)$ for all $j=1,...,J$. 
\end{defn}

Elements of $\mathcal{X}$ will be called \textit{patches}. Patches that are part of more than one budget plane will be called \textit{intersection patches}. Each budget can be uniquely expressed as union of patches; the number of patches that jointly comprise budget $\mathcal{B}_j$ will be called $I_j$. Note that $\sum_{j=1}^{J}I_{j} \geq I$, strictly so (because of multiple counting of intersection patches) if any two budget planes intersect.

\begin{rem}
\label{rem:i}
$I_j \leq I \leq 3^J$, hence $I_j$ and $I$ are finite.
\end{rem}

The partition $\mathcal{X}$ is the finite universal choice set alluded to earlier. The basic idea is that all choices from a given budget that are on the same patch induce the same directly revealed preferences, so are equivalent for the purpose of our test. Conversely, stochastic rationalizability does not at all constrain the distribution of demand on any patch. Therefore, rationalizability of $(P_1,\dots,P_J)$ can be decided by only considering the cross-sectional probabilities of patches on the respective budgets. We formalize this as follows. 

\begin{defn}
\label{defn:vector}
The \textit{vector representation} of $(\mathcal{B}_1,...,\mathcal{B}_J)$ is a $\left(\sum_{j=1}^{J}I_{j}\right)$-vector
$$ (x_{1|1},\dots,x_{I_1|1},x_{1|2},\dots,x_{I_J|J}),$$
where $(x_{1|j},\dots,x_{I_j|j})$ lists all patches comprising $\mathcal{B}_j$. The ordering of patches on budgets is arbitrary but henceforth fixed. Note that intersection patches appear once for each budget containing them.
\end{defn}
\begin{defn}
The \textit{vector representation} of $(P_1,...,P_J)$ is the $\left(\sum_{j=1}^{J}I_{j}\right)$-vector
$$ \pi:= (\pi_{1|1},\dots,\pi_{I_1|1},\pi_{1|2},\dots,\pi_{I_J|J}),$$
where $\pi_{i|j}:=P_j(x_{i|j})$.
\end{defn}

Thus, the vector representation of a stochastic demand system lists the probability masses that it assigns to patches.
\medskip

\textbf{Example \ref{example 1} continued.}
This example has a total of $5$ patches, namely the $4$ line segments identified in the figure and the intersection. The vector representations of $(\mathcal{B}_1,\mathcal{B}_2)$ and $(P_1,P_2)$ have $6$ components because the intersection patch is counted twice. If one disregards intersection patches (as we will do later), the vector representation of $(P_1,P_2)$ is $(\pi_{1|1},\pi_{2|1},\pi_{1|2},\pi_{2|2})$; see Figure \ref{figure 1}.  
\medskip

Next, a stochastic demand system is rationalizable iff it is a mixture of rationalizable nonstochastic demand systems. To intuit this, one may literally think of the latter as characterizing rational individuals. It follows that the vector representation of a rationalizable stochastic demand system must be the corresponding mixture of vector representations of rationalizable nonstochastic demand systems. Thus, define:

\begin{defn}
The \textit{rational demand matrix} $A$ is the (unique, up to ordering of columns) matrix such that the vector representation of each rationalizable nonstochastic demand system is exactly one column of $A$. The number of columns of $A$ is denoted $H$.
\end{defn}

\begin{rem}
\label{rem:h}
$H \leq \prod_{j=1}^J{I_j}$, hence $H$ is finite.
\end{rem}

We then have:

\begin{thm}
\label{prop:cone}
The following statements are equivalent:

\textbf{(i)} The stochastic demand system $(P_1,...,P_J)$ is rationalizable.

\textbf{(ii)} Its vector representation $\pi$ fulfills $\pi =A\nu$ for some $\nu \in \Delta^{H-1}$, the unit simplex in $\mathbf{R}^{H}$.

\textbf{(iii)} Its vector representation $\pi$ fulfills $\pi=A\nu$ for some $\nu \geq 0$.
\end{thm}

Theorem \ref{prop:cone} reduces the problem of testing \eqref{eq:model} to testing a null hypothesis about finite (though possible rather long) vector of probabilities. Furthermore, this hypothesis can be expressed as finite cone, a simple but novel observation that will be crucial for testing.\footnote{The idea of patches, as well as equivalence of (i) and (ii) in Theorem \ref{prop:cone}, were anticipated by \citeasnoun{mcfadden-2005}. While the explanation of patches is arguably unclear and (i)$\Leftrightarrow $(ii) is not explicitly pointed out, the idea is unquestionably there. The observation that (ii)$\Leftrightarrow $(iii) (more importantly: the idea of using this for testing) is new.}

We conclude this subsection with a few remarks.

\textbf{Simplification if demand is continuous.} Intersection patches are of lower dimension than budget planes. Thus, if the distribution of demand is continuous, their probabilities must be zero, and they can be eliminated from $\mathcal{X}$. This may considerably simplify $A$. Also, each remaining patch belongs to exactly one budget plane, so that $\sum_{j=1}^{J}I_{j}=I$. We impose this simplification henceforth and in our empirical application, but none of our results depend on it.

\textbf{GARP vs SARP.} Rationalizability of nonstochastic demand systems can be defined using either GARP or SARP. SARP will define a smaller matrix $A$, but nothing else changes. However, columns that are consistent with GARP but not SARP must select at least three intersection patches, so that GARP and SARP define the same $A$ if $\mathcal{X}$ was simplified to reflect continuous demand.

\textbf{Generality.} At its heart, Theorem \ref{prop:cone} only uses that choice from finitely many budgets reveals finitely many distinct revealed preference relations. Thus, it applies to any setting with finitely many budgets, irrespective of budgets' shapes. For example, the result was applied to kinked budget sets in \citeasnoun{Manski14} and could be used to characterize rationalizable choice proportions over binary menus, i.e. the Linear Order Polytope. The result furthermore applies to the ``random utility" extension of any other revealed preference characterization that allows for discretization of choice space; see \citeasnoun{DKQS16} for an example.

\subsection{Examples}

We next illustrate with a few examples. For simplicity, we presume continuous demand and therefore disregard intersection patches.
\medskip

\textbf{Example \ref{example 1} continued.}
Dropping the intersection patch, we have $I=4$ patches. Index vector representations as in Figure \ref{figure 1}, then the only excluded behavior is $(1,0,1,0)'$, thus
\begin{equation*}
A=\left( 
\begin{array}{ccc}
1 & 0 & 0 \\ 
0 & 1 & 1 \\ 
0 & 1 & 0 \\ 
1 & 0 & 1%
\end{array}
\right)
\begin{array}{r}
x_{1|1} \\
x_{2|1} \\
x_{1|2} \\
x_{2|2} \\
\end{array}.
\end{equation*}
The column cone of $A$ can be explicitly written as $\{(\nu_1,\nu_2+\nu_3,\nu_2,\nu_1+\nu_3)':\nu_1,\nu_2,\nu_3 \geq 0\}$. As expected, the only restriction on $\pi$ beyond adding-up constraints is that $\pi_{1|1}+\pi_{1|2}\leq 1$.

\begin{example} \label{example 2}
The following is the simplest example in which WARP does not imply SARP, so that applying Example \ref{example 1} to all pairs of budgets will only test a necessary condition. More subtly, it can be shown that the conditions in \citeasnoun{Kawaguchi} are only necessary as well. Let $K=J=3$ and assume a maximal pattern of intersection of budgets; for example, prices could be $( p_1,p_2,p_3) =((1/2,1/4,1/4),(1/4,1/2,1/4),(1/4,1/4,1/2))$. Each budget has $4$ patches for a total of $I=12$ patches, and one can compute
\begin{equation*}
\label{eq:ex2-A}
A=
\left(
\begin{array}{rrrrrrrrrrrrrrrrrrrrrrrrr}
0 & 0 & 0 & 0 & 1 & 0 & 0 & 0 & 0 & 0 & 0 & 0 & 0 & 0 & 0 & 1 & 0 & 0 & 0 & 0 & 0 & 1 & 0 & 0 & 0 \\ 
0 & 0 & 0 & 0 & 0 & 1 & 0 & 0 & 0 & 0 & 0 & 0 & 0 & 1 & 0 & 0 & 1 & 0 & 0 & 1 & 0 & 0 & 1 & 0 & 0 \\ 
0 & 1 & 0 & 0 & 0 & 0 & 1 & 0 & 0 & 0 & 0 & 1 & 0 & 0 & 0 & 0 & 0 & 1 & 0 & 0 & 0 & 0 & 0 & 1 & 0 \\ 
1 & 0 & 1 & 1 & 0 & 0 & 0 & 1 & 1 & 1 & 1 & 0 & 1 & 0 & 1 & 0 & 0 & 0 & 1 & 0 & 1 & 0 & 0 & 0 & 1 \\ 
0 & 0 & 0 & 0 & 0 & 0 & 0 & 0 & 1 & 0 & 0 & 0 & 0 & 1 & 1 & 0 & 0 & 0 & 0 & 0 & 0 & 0 & 0 & 0 & 0 \\ 
0 & 0 & 0 & 0 & 0 & 0 & 0 & 0 & 0 & 1 & 0 & 0 & 0 & 0 & 0 & 1 & 1 & 1 & 1 & 0 & 0 & 0 & 0 & 0 & 0 \\ 
1 & 0 & 0 & 1 & 0 & 0 & 0 & 0 & 0 & 0 & 1 & 0 & 0 & 0 & 0 & 0 & 0 & 0 & 0 & 1 & 1 & 0 & 0 & 0 & 0 \\ 
0 & 1 & 1 & 0 & 1 & 1 & 1 & 1 & 0 & 0 & 0 & 1 & 1 & 0 & 0 & 0 & 0 & 0 & 0 & 0 & 0 & 1 & 1 & 1 & 1 \\ 
1 & 1 & 1 & 0 & 0 & 0 & 0 & 0 & 0 & 0 & 0 & 0 & 0 & 0 & 0 & 0 & 0 & 0 & 0 & 0 & 0 & 0 & 0 & 0 & 0 \\ 
0 & 0 & 0 & 1 & 1 & 1 & 1 & 1 & 0 & 0 & 0 & 0 & 0 & 0 & 0 & 0 & 0 & 0 & 0 & 0 & 0 & 0 & 0 & 0 & 0 \\ 
0 & 0 & 0 & 0 & 0 & 0 & 0 & 0 & 1 & 1 & 1 & 1 & 1 & 0 & 0 & 0 & 0 & 0 & 0 & 0 & 0 & 0 & 0 & 0 & 0 \\ 
0 & 0 & 0 & 0 & 0 & 0 & 0 & 0 & 0 & 0 & 0 & 0 & 0 & 1 & 1 & 1 & 1 & 1 & 1 & 1 & 1 & 1 & 1 & 1 & 1
\end{array}
\right)
\begin{array}{r}
x_{1|1} \\ 
x_{2|1} \\ 
x_{3|1} \\ 
x_{4|1} \\ 
x_{1|2} \\ 
x_{2|2} \\ 
x_{3|2} \\ 
x_{4|2} \\ 
x_{1|3} \\ 
x_{2|3} \\ 
x_{3|3} \\ 
x_{4|3} 
\end{array} .
\end{equation*}

Interpreting $A$ requires knowing the geometry of patches. Given the ordering of patches used in the above, choice of $x_{3|1}$, $x_{2|2}$, and $x_{3|3}$ from their respective budgets would reveal a preference cycle, thus $A$ does not contain the column $(0,0,1,0,0,1,0,0,0,0,1,0)'$. We revisit this example in our Monte Carlo study.
\end{example}

\begin{figure}
\centering
\includegraphics[scale=1]{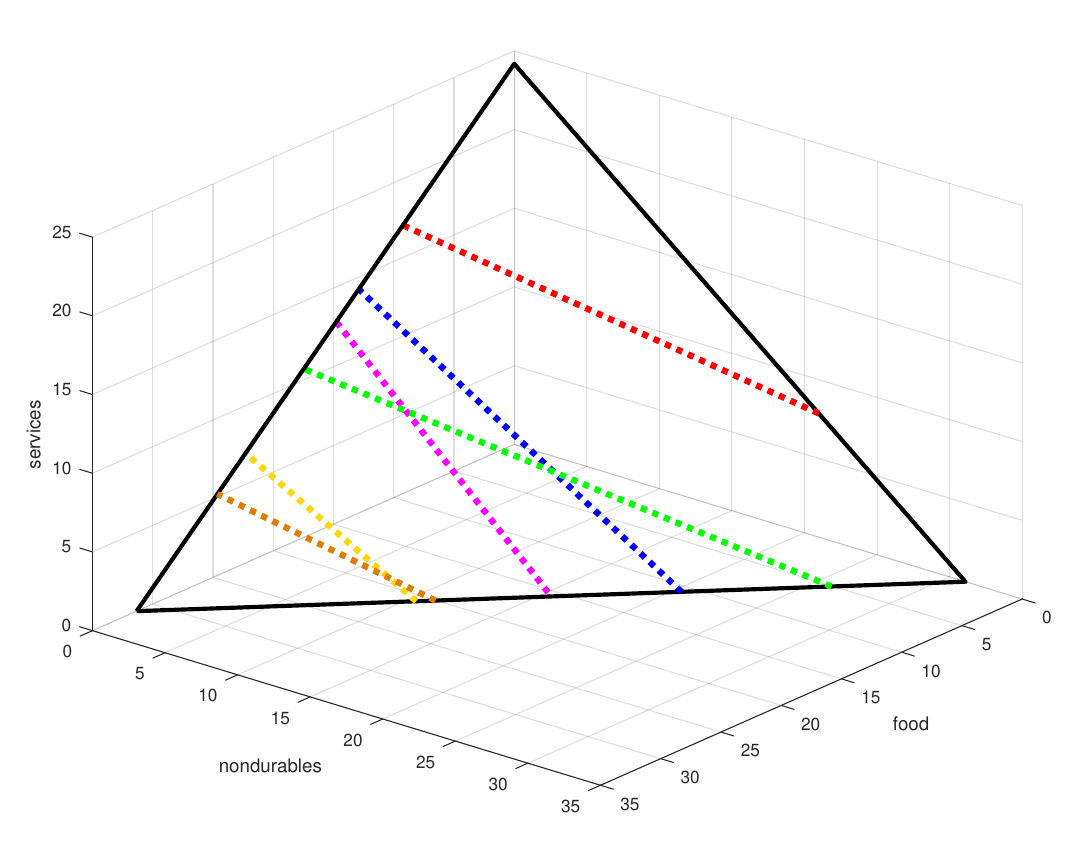}
\caption{Visualization of one budget set in the empirical application.}
\label{figure 2}
\end{figure}  
\begin{example}
Our empirical application has sequences of $J=7$ budgets in $\mathbf{R}^K$ for $K=3,4,5$ and sequences of $J=8$ budgets in $\mathbf{R}^3$. Figure \ref{figure 2} visualizes one of the budgets in $\mathbf{R}^3$ and its intersection with $6$ other budgets. There are total of $10$ patches (plus $15$ intersection patches). The largest $A$-matrices in the application are of sizes $78\times 336467$ and $79\times 313440$. In exploratory work using more complex examples, we also computed a matrix with over $2$ million columns.
\end{example}

\subsection{Computing A}
It should be clear now (and we formally show below) that the size of $A$, hence the cost of computing it, may escalate rapidly as examples get more complicated. We next elaborate how to compute $A$ from a vector of prices $(p_{1},...,p_{J})$. For ease of exposition, we drop intersection patches (thus SARP=GARP) and add remarks on generalization along the way. We split the problem into two subproblems, namely checking whether a binary ``candidate" vector $a$ is in fact a column of $A$ and finding all such vectors. 

\subsubsection{Checking rationalizability of $a$.} \label{sec:checkrat}

Consider any binary $I$-vector $a$ with one entry of $1$ on each subvector corresponding to one budget. This vector corresponds to a nonstochastic demand system. It is a column of $A$ if this demand system respects SARP, in which case we call $a$ \textit{rationalizable}.

To check such rationalizability, we initially extract a direct revealed preference relation over budgets. Specifically, if (an element of) $x_{i|j}$ is chosen from budget $\mathcal{B}_j$, then all budgets that are above $x_{i|j}$ are direct revealed preferred to $\mathcal{B}_j$. This information can be extracted extremely quickly.\footnote{\label{footnote-X}In practice, we compute a $(I \times J)$-matrix $X$ where, for example, the $i$-th row of $X$ is $(0,-1,1,1,1)$ if $x_{i|1}$ is on budget $\mathcal{B}_1$, below budget $\mathcal{B}_2$, and above the remaining budgets. This allows to vectorize construction of direct revealed preference relations, including strict vs. weak revealed preference, though we do not use the distinction.} 

We next exploit a well-known representation: Preference relations over $J$ budgets can be identified with directed graphs on $J$ labeled nodes by equating a directed link from node $i$ to node $j$ with revealed preference for $\mathcal{B}_i$ over $\mathcal{B}_j$. A preference relation fulfills SARP iff this graph is acyclic. This can be tested in quadratic time (in $J$) through a depth-first or breadth-first search. Alternatively, the Floyd-Warshall algorithm \cite{Floyd62} is theoretically slower but also computes rapidly in our application. Importantly, increasing $K$ does not directly increase the size of graphs checked in this step, though it allows for more intricate patterns of overlap between budgets and, therefore, potentially for richer revealed preference relations.

If intersection patches are retained, then one must distinguish between weak and strict revealed preference, and the above procedure tests SARP as opposed to GARP. To test GARP, one could use Floyd-Warshall or a recent algorithm that achieves quadatic time \cite{TSS15}.  

\subsubsection{Collecting rationalizable vectors.}
\label{sec:collect}

A total of $\prod_{j=1}^J I_J$ vectors $a$ could in principle be checked for rationalizability. Doing this by brute force rapidly becomes infeasible, including in our empirical application. However, these vectors can be usefully identified with the leaves (i.e. the terminal nodes) of a tree constructed as follows: (i) The root of the tree has no label and has $I_1$ children labelled $(x_{1|1},\dots,x_{I_1|1})$. (ii) Each child in turn has $I_2$ children labelled $(x_{1|2},\dots,x_{I_2|2})$, and so on for a total of $J$ generations beyond the root. Then there is a one-to-one mapping from leaves of the tree to conceivable vectors $a$, namely by identifying every leaf with the nonstochastic demand system that selects its ancestors. Furthermore, each non-terminal node of the tree can be identified with an incomplete $a$-vector that specifies choice only on the first $j<J$ budgets. The methods from Section \ref{sec:checkrat} can be used to check rationalizability of such incomplete vectors as well.

Our suggested algorithm for computing $A$ is a depth-first search of this tree. Importantly, rationalizability of the implied (possibly incomplete) vector $a$ is checked at each node that is visited. If this check fails, the node \textit{and all its descendants} are abandoned. A column of $A$ is discovered whenever a terminal node has been reached without detecting a choice cycle. Pseudocode for the tree search algorithm is displayed in Appendix B.

\subsubsection{Remarks on computational complexity.}

The cost of computing $A$ will escalate rapidly under any approach, but some meaningful comparison is possible. To do so, we consider three sequences, all indexed by $J$, whose first terms are displayed in Table \ref{table:sequences}. First, any two distinct nonstochastic demand systems induce distinct direct revealed preference relations; hence, $H$ is bounded above by $\bar H_J$, the number of distinct directed acyclic graphs on $J$ labeled nodes. This sequence -- and hence the worst-case cost of enumerating the columns of $A$, not to mention computing them -- is well understood \cite{Robinson73}, increases exponentially in $J$, and is displayed in the first row of Table \ref{table:sequences}.

Next, a worst-case bound on the number of terminal nodes of the aforementioned tree, hence on vectors that a brute force algorithm would check, is $2^{J(J-1)}$. This is simply because the number of conceivable candidate vectors $a$ equals $\prod_{j=1}^J I_j$, and every $I_j$ is bounded above by $2^{J-1}$. The corresponding sequence is displayed in the last row of Table \ref{table:sequences}.

Finally, some tedious combinatorial book-keeping (see Appendix C) reveals that the depth-first search algorithm visits at most $\sum_{j=2}^J \bar{H}_{j-1}2^{j(J+2-j)-2}$ nodes. This sequence is displayed in the middle row of the table, and the gain over a brute force approach is clear.\footnote{The comparison favors brute force because some nodes visited by a tree search are non-terminal, in which case rationalizability is easier to check. For example, this holds for $16$ of the $64$ nodes a tree search visits in Example \ref{example 2}.} 

It is easy to show that the ratio of any two sequences in the table grows exponentially. Also, all of the bounds are in principle attainable (though restricting $K$ may improve them) and are indeed attained in Examples \ref{example 1} and \ref{example 2}. In our empirical application, the bounds are far from binding (see Tables \ref{table:7 periods} and \ref{table:8 periods} for relevant values of $H$), but brute force was not always feasible, and the tree search improved on it by orders of magnitude in some cases where it was.

\begin{table}
\begin{tabular}{r|rrrrrrrr}
 J & 1 & 2 & 3 & 4 & 5 & 6 & 7 & 8 \\
 \hline
 $\bar{H}_J$ & $1$ & $3$ & $25$ & $543$ & $29281$ & $3.78 \times 10^6$ & $1.14 \times 10^9$ & $7.84 \times 10^{11}$ \\
 tree search & $1$ & $4$ & $64$ & $2048$ & $167936$ & $3.49 \times 10^7$ & $1.76 \times 10^{10}$ & $2.08 \times 10^{13}$ \\
 brute force & $1$ & $4$ & $64$ & $4096$ & $1048576$ & $1.08 \times 10^8$ & $4.40 \times 10^{12}$ & $7.21 \times 10^{16}$  \\
 \multicolumn{9}{c}{} \\
\end{tabular}
\caption{Worst-case bounds on the number of rationalizable vectors $a$ and on the number of candidate vectors visited by different algorithms.}
\label{table:sequences}
\end{table}

\subsubsection{Further Refinement.}

A modest amount of problem-specific adjustment may lead to further improvement. The key to this is contained in the following result.
\begin{thm}
\label{prop1}Suppose that not all budget planes mutually intersect; in particular, there exists $M<J$ s.t. $\mathcal{B}_{J}$ is either above all or below all of $(\mathcal{B}_{1},...,\mathcal{B}_{M})$. Suppose also that choices from $(\mathcal{B}_{1},...,\mathcal{B}_{J-1})$ are rationalizable. Then choices from $(\mathcal{B}_{1},...,\mathcal{B}_{J})$ are rationalizable iff choices from $(\mathcal{B}_{M+1},...,\mathcal{B}_{J})$ are.
\end{thm}
If the geometry of budgets allows it -- this is particularly likely if budgets move outward over time and even guaranteed if some budget planes are parallel -- Theorem \ref{prop1} can be used to construct columns of $A$ recursively from columns of $A$-matrices that correspond to a smaller $J$. The gain can be tremendous because, at least with regard to worst-case cost, one effectively moves one or more columns to the left in Table \ref{table:sequences}. A caveat is that application of Theorem \ref{prop1} may require to manually reorder budgets so that it applies. Also, while the internal ordering of $(\mathcal{B}_1,\dots,\mathcal{B}_M)$ and $(\mathcal{B}_{M+1},\dots,\mathcal{B}_{J-1})$ does not matter, the theorem may apply to distinct partitions of the same set of budgets. In that case, any choice of partition will accelerate computations, but we have no general advice on which is best. We tried the refinement in our empirical application, and it considerably improved computation time for some of the largest matrices. However, the tree search proved so fast that, in order to keep it transparent, our replication code omits this step.

\section{Statistical Testing}

\label{sec:Inference}

This section lays out our statistical testing procedure in the idealized situation where, for finite $J$, repeated cross-sectional observations of demand over $J$ periods are available to the econometrician. Formally, for each $1 \leq j \leq J$, suppose one observes $N_j$ random draws of $y$ distributed according to $P_j$ defined in \eqref{eq:define_pi_j}. Define $N = \sum_{j=1}^J N_j$ for later use. Clearly, $P_j$ can be estimated consistently as $N_j \uparrow \infty$ for each $j$, $1 \leq j \leq J$. The question is whether the estimated distributions may, up to sampling uncertainty, have arisen from a RUM. We define a test statistic and critical value and show that the resulting test is uniformly asymptotically valid over an interesting range of d.g.p.'s.

\subsection{Null Hypothesis and Test Statistic}
By Theorem \ref{prop:cone}, we wish to test:

\ 

\noindent (\textbf{H}$_A$): There exist $\nu \geq 0$ such that $A\nu =\pi$.

\ 

\noindent This hypothesis is equivalent to

\

\noindent (\textbf{H}$_{B}$): \quad $\min_{\eta \in \mathcal C}[\pi -\eta ]^{\prime
}\Omega \lbrack \pi -\eta ]=0$,

\ 

\noindent where $\Omega $ is a positive definite matrix (restricted to be
diagonal in our inference procedure) and $\mathcal{C}:=\{A\nu |\nu \geq 0\}$ is a convex cone in $\mathbf{R}^{I}$. The solution $%
\eta _{0}$ of (\textbf{H}$_{B}$) is the projection of $\pi \in \mathbf{R}%
_{+}^{I}$ onto $\mathcal{C}$ under the weighted norm $\Vert x\Vert _{\Omega }=\sqrt{%
	x^{\prime }\Omega x}$. The corresponding value of the objective function is
the squared length of the projection residual vector. The projection $\eta
_{0}$ is unique, but the corresponding $\nu $ is not. Stochastic rationality
holds if and only if the length of the residual vector is zero.

A natural sample counterpart of the objective function in (\textbf{H}$_{B}$)
would be $\min_{\eta \in \mathcal C}[\hat{\pi}-\eta ]^{\prime }\Omega \lbrack \hat{\pi}-\eta ]$, where $\hat{\pi}$ estimates $\pi $, for example by sample choice
frequencies. The usual scaling yields
\begin{eqnarray}
	\mathcal{J}_N &:=&N\min_{\eta \in \mathcal C}[\hat{\pi}-\eta ]^{\prime }\Omega \lbrack \hat{%
		\pi}-\eta ]  \label{J-def} \\
	&=&N\min_{\nu \in \mathbf{R}_{+}^{H}}[\hat{\pi}-A\nu ]^{\prime }\Omega
	\lbrack \hat{\pi}-A\nu ].  \notag
\end{eqnarray}
Once again, $\nu $ is not unique at the optimum, but $\eta =A\nu $ is. Call its optimal value $\hat{\eta}$. Then $\hat{\eta}=\hat{\pi}$, and $\mathcal{J}_N=0$, if the estimated choice probabilities $\hat{\pi}$ are stochastically rationalizable; obviously, our null hypothesis will be accepted in this case.

\subsection{Simulating a Critical Value.}
\label{sec:bootstrap}
We next explain how to get a valid critical value for $\mathcal{J}_N$ under the assumption that $\hat{\pi}$ estimates the probabilities of patches by corresponding sample frequencies and that one has $R$ bootstrap replications $\hat{\pi}^{\ast (r)},r=1,...,R$. Thus, $\hat{\pi}^{\ast (r)}-\hat{\pi}$ is a natural bootstrap analog of $\hat{\pi}-\pi$. We will make enough assumption to ensure that its distribution consistently estimates the distribution of $\hat{\pi}-\pi_0$, where $\pi_0$ is the true value of $\pi$. The main difficulty is that one cannot use $\hat{\pi}$ as bootstrap analog of $\pi_0$.

Our bootstrap procedure relies on a tuning parameter $\tau_N$ chosen s.t. $\tau _{N}\downarrow 0$ and $\sqrt{N}%
\tau _{N}\uparrow \infty $.\footnote{In this section's simplified setting and if $\hat{\pi}$ collects sample frequencies, a reasonable choice would be 
	\begin{equation*}
		\tau _{N}=\sqrt{\frac{\log \underline{N}}{\underline{N}}}
	\end{equation*}%
	where $\underline{N}=\min_{j}N_{j}$ and $N_{j}$ is the number of observations on Budget $\mathcal{B}_{j}$: see \eqref{eq:pihat}.  This choice corresponds to the ``BIC choice'' in \citeasnoun{andrews-soares}.  We will later propose a different $\tau_N$ based on how $\pi$ is in fact estimated.} Also, we restrict $\Omega$ to be diagonal and positive definite and let $\mathbf{1}_{H}$ be a $H$-vector of ones\footnote{In principle, $\mathbf{1}_{H}$ could be any strictly positive $H$-vector, though a data based choice of such a vector is beyond the scope of the paper.}.   The restriction on $\Omega$ is important: Together with a geometric feature of the column vectors of the matrix $A$, it ensures that constraints which are fulfilled but with small slack become binding through the \textit{Cone Tightening} algorithm we are about to describe. A non-diagonal weighting matrix can disrupt this property. For further details on this point and its proof, the reader is referred to Appendix A. Our
procedure is as follows:

\begin{itemize}
	\item[(i)] Obtain the \textit{$\tau_N$-tightened restricted} estimator $\hat{%
		\eta}_{\tau _{n}}$, which solves%
	\begin{eqnarray*}
	\min_{\eta \in \mathcal{C}_{\tau _{N}}}N[\hat{\pi}-\eta ]^{\prime
	}\Omega \lbrack \hat{\pi}-\eta ] =\min_{[\nu -\tau _{N}\mathbf{1}_{H}/H]\in \mathbf{R}_{+}^{H}}N[\hat{\pi}%
		-A\nu ]^{\prime }\Omega \lbrack \hat{\pi}-A\nu ].
	\end{eqnarray*}
	
	\item[(ii)] Define the \textit{$\tau_N$-tightened recentered} bootstrap estimators
	\begin{equation*}
		\hat \pi^{*(r)}_{\tau_N} := \hat \pi^{*(r)} - \hat \pi + \hat \eta_{\tau_N},
		\quad r = 1,...,R.
	\end{equation*}
	
	\item[(iii)] The bootstrap test statistic is
	\begin{eqnarray*}
		\mathcal{J}_N^{*(r)} = \min_{[\nu - \tau_N \mathbf{1}_H/H] \in \mathbf{R}_+^H}N[\hat
		\pi^{*(r)}_{\tau_N} - A \nu]^{\prime }\Omega[\hat \pi^{*(r)}_{\tau_N} - A \nu%
		],
	\end{eqnarray*}
	for $r = 1,...,R$.
	
	\item[(iv)] Use the empirical distribution of $\mathcal{J}_N^{\ast (r)}, r=1,...,R$ to obtain the critical value for $\mathcal{J}_N$.
\end{itemize}

\ 

The object $\hat \eta_{\tau_N}$ is the true value of $\pi$ in the bootstrap population, i.e. it is the bootstrap analog of $\pi_0$. It differs from $\hat{\pi}$ through a ``double recentering.'' To disentangle the two recenterings, suppose first that $\tau_N=0$. Then inspection of step (i) of the algorithm shows that $\hat{\pi}$ would be projected onto the cone $\mathcal{C}$. This is a relatively standard recentering ``onto the null'' that resembles recentering of the $J$-statistic in overidentified GMM. However, with $\tau_N>0$, there is a second recentering because the cone $\mathcal{C}$ itself has been tightened. We next discuss why this recentering is needed.

\subsection{Discussion}
Our testing problem is related to the large literature on inequality testing but adds an important twist. Writing $\{a_{1},a_{2},...,a_{H}\}$ for the
column vectors of $A$, one has 
\begin{equation*}
	\mathcal{C}=\mathrm{{cone}(A)}:=\{\nu _{1}a_{1}+...+\nu _{H}a_{H}:\nu _{h}\geq 0\},
\end{equation*}
i.e. the set $\mathcal{C}$ is a finitely generated cone.  The following result, known as the {\textsc{Weyl-Minkowski Theorem}}, provides an alternative representation that is useful for theoretical developments of our statistical testing procedure.\footnote{See \citeasnoun{Gruber}, \citeasnoun{Grunbaum03}, and \citeasnoun{Ziegler}, especially Theorem 1.3, for these results and other 	materials concerning convex polytopes used in this paper.}  
\begin{thm} \label{thm:WM} (Weyl-Minkowski Theorem for Cones)  A subset $\mathcal{C}$ of ${\bf R}^I$
	is a finitely generated cone	 
	\begin{equation}\label{eq:v-rep} 
		\mathcal{C}=\{\nu _{1}a_{1}+...+\nu _{H}a_{H}:\nu _{h}\geq 0\} \text{ for some } A = [a_1,...,a_H] \in {\bf R}^{I\times H}
	\end{equation}
	if and only if it is a finite intersection of closed half spaces 
	\begin{equation}\label{eq:h-rep}
		\mathcal{C}=\{t\in \mathbf{R}^{I}|Bt\leq 0\} 	\text{ for some } B \in {\bf R}^{m\times I}.
	\end{equation}
	The expressions in \eqref{eq:v-rep} and \eqref{eq:h-rep} are called a $\mathcal{V}$-representation (as in ``vertices'') and a $\mathcal{H}$-representation (as in ``half spaces'') of $\mathcal{C}$, respectively.
\end{thm}
The ``only if'' part of the theorem (which is \textsc{Weyl's Theorem}) shows that our rationality hypothesis
$
\pi \in \mathcal{C}, \mathcal{C} = \{A\nu|\nu \geq 0\}
$
in terms of a $\mathcal{V}$-representation can be re-formulated in an $\mathcal{H}$-representation using an appropriate matrix $B$, at least in theory.  If such $B$ were available, our testing problem would resemble tests of%
\begin{equation*}
	H_{0}:B\theta \geq 0\quad B\in \mathbf{R}^{p\times q}\text{ is known}
\end{equation*}%
based on a quadratic form of the empirical discrepancy between $B\theta$ and $\eta$ minimized over $\eta \in {\bf R}_+^q$.   
This type of problem has been studied extensively; see references in Section \ref{sec:literature}. Its analysis is intricate because the limiting distribution of such a statistic depends discontinuously on the true value of $B\theta$. One common way to get a critical value is to consider the globally least favorable case, which is $\theta =0$. A less conservative strategy widely followed in the econometric literature on moment inequalities is \textit{Generalized Moment Selection} (GMS; see \citeasnoun{andrews-soares}, \citeasnoun{bugni-2010}, \citeasnoun{Canay10}). If we had the $\mathcal{H}$-representation of $\mathcal{C}$, we might conceivably use the same technique. However, the duality between the two representations is purely theoretical: In practice, $B$ cannot be computed from $A$ in high-dimensional cases like our empirical application.

We therefore propose a \textit{tightening} of the cone $\mathcal{C}$ that is computationally feasible and will have a similar effect as GMS. The idea is to tighten the constraint on $\nu$ in \eqref{J-def}. In particular, define $\mathcal{C}_{\tau _{N}}:=\{A\nu |\nu \geq \tau _{N}\mathbf{1}_{H}/H\}$ and define $\hat{\eta}_{\tau _{N}}$ as 
\begin{eqnarray}
	\hat{\eta}_{\tau _{N}} &:=&\argmin_{\eta \in \mathcal{C}_{\tau _{N}}}N[\hat{\pi}-\eta ]^{\prime
	}\Omega \lbrack \hat{\pi}-\eta ] \\
	&=&\argmin_{A\nu: \: [\nu -\tau _{N}\mathbf{1}_{H}/H]\in \mathbf{R}_{+}^{H}}N[\hat{\pi}%
	-A\nu ]^{\prime }\Omega \lbrack \hat{\pi}-A\nu ].
	\nonumber
\end{eqnarray}%
Our proof establishes that constraints in the $\mathcal{H}$-representation that are almost binding at the original problem's solution
(i.e., their slack is difficult to be distinguished from zero at the sample
size) will be binding with zero slack after tightening. Suppose that $\sqrt{N}(\hat{\pi}-\pi )\rightarrow _{d}N(0,S)$ and let $\hat{S}$ consistently estimate $S$. Let $\tilde{\eta}_{\tau _{N}}:=\hat{\eta}_{\tau _{N}}+\frac{1}{\sqrt{N}}N(0,\hat{S})$ or a bootstrap random variable and use the distribution of 
\begin{eqnarray}
	\tilde{\mathcal J}_{N} &:=&\min_{\eta \in \mathcal{C}_{\tau _{N}}}N[\tilde{\eta}%
	_{\tau _{N}}-\eta ]^{\prime }\Omega \lbrack \tilde{\eta}_{\tau _{N}}-\eta ]
	\label{J-tilde} \\
	&=&\min_{[\nu -\tau _{N}\mathbf{1}_{H}/H]\in \mathbf{R}_{+}^{H}}N[\tilde{\eta}%
	_{\tau _{N}}-A\nu ]^{\prime }\Omega \lbrack \tilde{\eta}_{\tau _{N}}-A\nu ],
	\notag
\end{eqnarray}%
to approximate the distribution of $\mathcal{J}_N$. This has the same theoretical justification as the inequality selection procedure. Unlike the latter, however, it avoids the use of an $\mathcal{H}$-representation, thus offering a computationally feasible testing procedure.

To further illustrate the duality between $\mathcal{H}$- and $\mathcal{V}$-representations, we revisit the first two examples. It is not possible to compute $B$-matrices in our empirical application.
\smallskip

\textbf{Example \ref{example 1} continued.} With two intersecting budget planes, the cone $\mathcal{C}$ is represented by 
\begin{equation}
B=~
\begin{blockarray}{rrrr}
x_{1|1} & x_{2|1} & x_{1|2} & x_{2|2}~ \\
\begin{block}{(rrrr)}
-1 & 0 & 0 & 0~ \\ 
0 & 0 & -1 & 0~ \\ 
-1 & -1 & 1 & 1~ \\ 
1 & 1 & -1 & -1~ \\ 
1 & 0 & 0 & -1~ \\
\end{block}
\end{blockarray}~.
\label{eq:b-matrix}
\end{equation}
The first two rows of $B$ are nonnegativity constraints (the other two such constraints are redundant), the next two rows are an equality constraint forcing the sum of probabilities to be constant across budgets, and only the last constraint is a substantive economic constraint. If the estimator $\hat{\pi}$ fulfills the first four constraints by construction, then the testing problem simplifies to a test of $(1,0,0,-1)\pi \leq 0$, the same condition identified earlier.
\smallskip

\textbf{Example \ref{example 2} continued.} Eliminating nonnegativity and adding-up constraints for brevity, numerical evaluation reveals
\begin{equation}
\label{eq:ex2-matrix}
B=~
\begin{blockarray}{rrrrrrrrrrrr}
x_{1|1} & x_{2|1} & x_{3|1} & x_{4|1} & x_{1|1} & x_{1|2} & x_{3|1} & x_{4|1} & x_{1|1} & x_{2|1} & x_{3|1} & x_{4|1}~ \\
\begin{block}{(rrrrrrrrrrrr)}
1 & 1 & 0 & 0 & 0 & -1 & 0 & -1 & 0 & 0 & 0 & 0~ \\ 
0 & -1 & 0 & -1 & 0 & 0 & 0 & 0 & 1 & 1 & 0 & 0~ \\ 
0 & 0 & 0 & 0 & 1 & 1 & 0 & 0 & 0 & -1 & 0 & -1~ \\
0 & -1 & -1 & -1 & 1 & 0 & 0 & 0 & 1 & 0 & 0 & 0~ \\
1 & 1 & 1 & 0 & 0 & 0 & 0 & -1 & 0 & 0 & 0 & -1~ \\
\end{block}
\end{blockarray}~. 
\end{equation}
The first three rows are constraints on pairs of budgets that mirror the last row of \eqref{eq:b-matrix}. The next two constraints are not implied by these, nor by additional constraints in \citeasnoun{Kawaguchi}, but they imply the latter.

\subsection{Theoretical Justification}
We now provide a detailed justification. First, we formalize the notion that choice probabilities are estimated by sample frequencies. Thus, for
each budget set $\mathcal{B}_j$, denote the choices of $N_{j}$ individuals, indexed by $%
n=1,...,N_{j}$, by 
\begin{equation*}
	d_{i|j,n}=%
	\begin{cases}
		& \mbox{$1$ if individual $n$ chooses $x_{i|j}$} \\ 
		& \mbox{$0$ otherwise}%
	\end{cases}%
	\qquad n=1,...,N_{J}.
\end{equation*}%
Assume that one observes $J$ random samples $\{\{d_{i|j,n}\}_{i=1}^{I_{j}}%
\}_{n=1}^{N_{j}}$, $j=1,2,...,J$. For later use, define 
\begin{equation*}
	d_{j,n}:=\left[ 
	\begin{array}{c}
		d_{1|j,n} \\ 
		\vdots \\ 
		d_{I_{j}|j,n}%
	\end{array}%
	\right] ,\quad N=\sum_{j=1}^{J}N_{J}.
\end{equation*}%
An obvious way to estimate the vector $\pi $ is to use choice frequencies 
\begin{equation}
	\hat{\pi}_{i|j}=\sum_{n=1}^{N_{j}}d_{i|j,n}/N_{j},i=1,...,I_{j},j=1,...,J.
	\label{eq:pihat}
\end{equation}%

The next lemma, among other things, shows that our tightening of the $\mathcal{V}$-representation of $\mathcal{C}$ is equivalent to a tightening its $\mathcal{H}$-representation but leaving $B$ unchanged. For a matrix $B$, let $\mathrm{col}(B) $ denote its column space.

\begin{lem}
	\label{lem:c_tau} For $A\in \mathbf{R}^{I\times H}$, let
	$$
	\mathcal{C}  = \{A \nu|\nu \geq 0\}.
	$$ 
	Also let 
	\begin{equation*}
		\mathcal{C}=\{t:Bt\leq 0\}
	\end{equation*}%
	be its $\mathcal H$-representation for some $B \in \mathbf{R}^{m\times I}$ such that  
	$B = \begin{bmatrix}
	B^\leq \\ B^=
	\end{bmatrix}$, where the submatrices $B^\leq \in \mathbf{R}^{{\bar m} \times I}$ and $B^= \in \mathbf{R}^{(m - \bar m) \times I}$ correspond to inequality and equality constraints, respectively.  For $\tau >0$ define 
	$$
	\mathcal{C}_\tau  = \{A \nu|\nu \geq (\tau/H) \mathbf{1}_{H}\}.
	$$ 
	Then one also has 
	\begin{equation*}
		\mathcal{C}_\tau=\{t:Bt\leq -\tau \phi \}
	\end{equation*}%
	for some $\phi=(\phi_1,...,\phi_m)' \in \mathrm{col}(B)$ with the properties that (i) $\bar \phi :=  	[\phi_1,...,\phi_{\bar m}]'  \in {\bf R}^{\bar m}_{++}$, and (ii) $\phi_k = 0$ for $k > \bar m$.
\end{lem}
Lemma \ref{lem:c_tau} is {\it not} just a re-statement of the {\textsc{Minkowski-Weyl theorem}} for polyhedra, which would simply say $\mathcal{C}_\tau  = \{A \nu|\nu \geq (\tau/H) \mathbf{1}_{H}\}$ is alternatively represented as an intersection of closed halfspaces.  The lemma instead shows that the inequalities in the $\mathcal H$-representation becomes tighter by $\tau \phi$ after tightening the $\mathcal V$-representation by $\tau_N {\bf 1}_H/H$, with the same matrix of coefficients $B$ appearing both for $\mathcal{C}$ and $\mathcal{C}_\tau$. Note that for notational convenience, we rearrange rows of $B$ so that the genuine inequalities come first and pairs of inequalities that represent equality constraints come last.\footnote{In the matrix displayed in \eqref{eq:b-matrix}, the third and fourth row would then come last.} This is w.l.o.g.; in particular, the researcher does not need to know which rows of $B$ these are.  Then as we show in the proof, the elements in $\phi$ corresponding to the equality constraints are automatically zero when we tighten the space for {\it all} the elements of $\nu$ in the $\mathcal V$-representation.  This is a useful feature that makes our methodology work in the presence of equality constraints.    

The following assumptions are used for our asymptotic theory.

\begin{ass}
	\label{ass1} For all $j = 1,...,J$, $\frac {N_j}{N} \rightarrow \rho_j$ as $%
	N \rightarrow \infty$, where $\rho_j >0$.
\end{ass}

\begin{ass}\label{ass:samping}
	$J$ repeated cross-sections of random samples $\left\{ 
	\{d_{i|j,n(j)}\}_{i=1}^{I_{j}} \right\} _{n(j)=1}^{N_{j}},j=1,...,J$, are observed. 
\end{ass}
\noindent The econometrician also observes the normalized price vector $p_j$, which is fixed in this section, for each $1 \leq j \leq J$.  Let ${\mathcal{P}}$ denote the set of all 
$\pi$'s that satisfy Condition \ref{condition 1} in Appendix A for some (common) value of $%
(c_1,c_2)$.

\begin{thm}
	\label{thm1} Choose $\tau _{N}$ so that $\tau _{N}\downarrow 0$ and $\sqrt{N}%
	\tau _{N}\uparrow \infty $. Also, let $\Omega $ be diagonal, where all the
	diagonal elements are positive. Then under Assumptions \ref{ass1} and \ref{ass:samping} 
	\begin{equation*}
		\liminf_{N\rightarrow \infty }\inf_{\pi \in \mathcal{P}\cap \mathcal{C}}\Pr
		\{\mathcal{J}_N\leq \hat{c}_{1-\alpha }\}=1-\alpha
	\end{equation*}%
	where $\hat{c}_{1-\alpha }$ is the $1-\alpha $ quantile of $\tilde{\mathcal J}_N$, $0\leq \alpha \leq \frac{1}{2}$.
\end{thm}
While it is obvious that our tightening contracts the cone, the result depends on a more delicate feature, namely that we (potentially) turn non-binding inequalities from the $\mathcal{H}$-representation into
binding ones but not vice versa. This feature is not universal to cones as they get contracted. Our proof establishes that it generally obtains if $\Omega $ is the identity
matrix and all corners of the cone are acute. In this paper's application,
we can further exploit the cone's geometry to extend the result to any
diagonal $\Omega $.\footnote{It is possible to replace $\Omega$ with its consistent estimator and retain uniform asymptotic validity, if we further impose a restriction on the class of distributions over which we define the size of our test.  Note, however, that our $\mathcal P$ in our Theorem \ref{thm1} (and its variants in Theorems \ref{thmSM} and \ref{thmEC}) allows for some elements of the vector $\pi$ being zeros.  This makes the use of the reciprocals of estimated variances for the diagonals of the weighting matrix potentially problematic, as it invalidates the asymptotic uniform validity since the required triangular CLT does not hold under parameter sequences where the elements of $\pi$ converge to zeros.  The use of fixed $\Omega$, which we recommend in implementing our procedure, makes contributions from these terms asymptotically negligible, thereby circumventing this problem.}   Our method immediately applies to other testing problems featuring $\mathcal{V}$%
-representations if analogous features can be verified.

\section{Methods for Typical Survey Data}
\label{sec:extending}

The methodology outlined in Section \ref{sec:Inference} requires that (i) the observations available to the econometrician are drawn on a finite number of budgets and (ii) the budgets are given exogenously, that is, unobserved heterogeneity and budgets are assumed to be independent.  These conditions are naturally satisfied in some applications.  The empirical setting in Section \ref{sec:empirics}, however, calls for modifications because Condition (i) is certainly violated in it and imposing Condition (ii) would be very restrictive.   These are typical issues for  a survey data set.  This section addresses them.

Let $P_u$ denote the marginal probability law of $u$, which we assume does not depend on $j$. We do not, however, assume that the laws of other random elements, such as income, are time homogeneous.  Let $w = \log(W)$ denote log total expenditure, and suppose the researcher chooses a value $\underline{w}_j$ for $w$ for each period $j$.  Note that our algorithm and asymptotic theory remain valid if multiple values of $w$ are chosen for each period.
Let $w_{n(j)}$ be
the log total expenditure of consumer $n(j)$, $1\leq n(j)\leq N_{j}$ observed in period $j$.
\begin{ass}\label{ass:sampingsmooth}
	$J$ repeated cross-sections of random samples $\left\{ \left(
	\{d_{i|j,n(j)}\}_{i=1}^{I_{j}},w_{n(j)}\right) \right\} _{n(j)=1}^{N_{j}},j=1,...,J$, are observed. 
\end{ass}
\noindent The econometrician also observes the unnormalized price vector $\tilde p_j$, which is fixed, for each $1 \leq j \leq J$.     

We first assume that the total expenditure is exogenous, in the sense that 
$
w \indep u
$
holds under every $P^{(j)}, 1 \leq j \leq J$.  This exogeneity assumption will be relaxed shortly.   Let
$\pi_{i|j}(w):=\Pr\{d_{i|j,n(j)}=1|w_{n(j)}=  \underline{w}_j \}$ and 
writing $\pi _{j}:=(\pi _{1|j},...,\pi _{I_{j}|j})^{\prime }$
and $\pi :=(\pi _{1}^{\prime },...,\pi _{J}^{\prime })^{\prime }=(\pi
_{1|1},\pi _{2|1},...,\pi _{I_{J}|J})^{\prime }$, the stochastic rationality condition is given by 
$
\pi \in \mathcal{C}
$
as before.  Note that this $\pi$ can be estimated by standard nonparametric procedures.  For concreteness, we use a series estimator, as defined and analyzed in Appendix A.  
The smoothed version of $%
\mathcal{J}_N$  (also denoted $\mathcal J_N$ for simplicity)  is obtained using the series estimator for $\hat{\pi}$ in %
\eqref{J-def}. 
In  Appendix A we also present an algorithm for obtaining the bootstrapped version $\tilde {\mathcal J}_N$ of the smoothed statistic.

In what follows, $F_{j}$ signifies the joint distribution of $%
(d_{i|j,n(j)},w_{n(j)})$. Let $\mathcal{F}$ be the set of all $%
(F_{1},...,F_{J}) $ that satisfy Condition \ref{condition:kernel}  in Appendix A for some $%
(c_1,c_2,\delta,\zeta(\cdot ))$.

\begin{thm}
	\label{thmSM}  Let Condition \ref{cond:tuning1} hold.  Also let $\Omega$ be diagonal where all the
	diagonal elements are positive. Then under Assumptions \ref{ass1} and \ref{ass:sampingsmooth} 
	\begin{equation*}
		\liminf_{N \rightarrow \infty} \inf_{(F_1,...,F_J) \in \mathcal{F}} \Pr\{\mathcal{J}_N
		\leq \hat c_{1 - \alpha}\} = 1 - \alpha
	\end{equation*}
	where $\hat c_{1 - \alpha}$ is the $1 - \alpha$ quantile of $\tilde{\mathcal{J}}_N$, $0 \leq \alpha \leq \frac 1 2$.
\end{thm}

Next, we  relax the assumption that consumer's utility functions are realized independently from $W$. 
For each fixed value $\underline{w}_j$ and the unnormalized price vector $\tilde p_j$ in period $j$,  $1 \leq j \leq J$, define the endogeneity corrected conditional probability{\footnote{This is the conditional choice probability if $p$ is (counterfactually) assumed to be exogenous. We call it ``endogeneity corrected" instead of ``counterfactual" to avoid confusion with rationality constrained, counterfactual prediction.}
	\begin{eqnarray*}
		\pi (\tilde{p}_j/e^{\underline{w}_j},x_{i|j}) &:=& 
		\int {\bf 1}\{D_j(\underline{w}_j,u) \in x_{i|j} \}  dP_u
	\end{eqnarray*}
	where $D_j(w,u) := D(\tilde{p}_j/e^w,u)$.
	Then Theorem \ref{prop:cone} still applies to 
	$$
	\pi_{\rm{EC}} := [\pi(p_1,x_{1|1}),...,\pi(p_{1},x_{I_1|1}),\pi(p_2,x_{1|2}),...,\pi(p_{2},x_{I_2|2}),...,\pi(p_J,x_{1|J}),...,\pi(p_{J},x_{I_J|J})]'.
	$$ 
	Suppose there exists a control variable $\varepsilon$ such that 
	$
	{w \indep u | \varepsilon}
	$
	holds under every $P^{(j)}, 1 \leq j \leq J$.  See \eqref{eq:cfexample} in Section A for an example.   We propose to use a fully nonparametric, control function-based two-step estimator, denoted by $\widehat{\pi_{\mathrm{EC}}}$, to define our endogeneity-corrected test statistic ${\mathcal J}_{\mathrm{EC}_N}$; see Appendix A for details.  For this, the bootstrap procedure needs to be adjusted appropriately to obtain the bootstrapped statistic $\tilde {\mathcal J}_{\mathrm{EC}_N}$: once again, the reader is referred to Appendix A.   This is the method we use for the empirical results reported in Section \ref{sec:Empirical Application}.   Let $z_{n(j)}$ be the $n(j)$-th observation of the instrumental variable $z$ in period $j$.  
\begin{ass}\label{ass:sampingendogeneity}
	$J$ repeated cross-sections of random samples $\left\{ \left(
	\{d_{i|j,n(j)}\}_{i=1}^{I_{j}},x_{n(j)},z_{n(j)}\right) \right\} _{n=1}^{N_{j}},j=1,...,J$, 
	 are observed. 
\end{ass}
\noindent The econometrician also observes the unnormalized price vector $\tilde p_j$, which is fixed, for each $1 \leq j \leq J$.

	In what follows, $F_{j}$ signifies the joint distribution of $%
	(d_{i|j,n(j)},w_{n(j)},z_{n(j)})$. Let $\mathcal{F}_{\rm EC}$ be the set of all $%
	(F_{1},...,F_{J}) $ that satisfy Condition \ref{condition:EC} in Appendix A for some $(c_1,c_2,\delta_1,\delta,\zeta_r(\cdot),\zeta_s(\cdot),\zeta_1(\cdot))$. Then we have:
	\begin{thm}
		\label{thmEC} Let Condition \ref{cond:tuning2} hold.    Also let $\Omega$ be diagonal where all the diagonal elements are positive. Then under Assumptions \ref{ass1} and \ref{ass:sampingendogeneity}, 
		\begin{equation*}
			\liminf_{N \rightarrow \infty} \inf_{(F_1,...,F_J) \in \mathcal{F}_{\rm EC}} \Pr\{\mathcal J_{{\rm EC}_N}
			\leq \hat c_{1 - \alpha}\} = 1 - \alpha
		\end{equation*}
		where $\hat c_{1 - \alpha}$ is the $1 - \alpha$ quantile of $\tilde{\mathcal
		J}_{{\rm EC}_N}$, $0 \leq \alpha \leq \frac 1 2$.
	\end{thm}

\section{Monte Carlo Simulations}
\label{sec:monte carlo}

We next analyze the performance of Cone Tightening in a small Monte Carlo study. To keep examples transparent and to focus on the core novelty, we model the idealized setting of Section \ref{sec:Inference}, i.e. sampling distributions are multinomial over patches. In addition, we focus on Example \ref{example 2}, for which an $\mathcal{H}$-representation in the sense of Weyl-Minkowski duality is available; see displays \eqref{eq:ex2-A} and \eqref{eq:ex2-matrix} for the relevant matrices.\footnote{This is also true of Example \ref{example 1}, but that example is too simple because the test reduces to a one-sided test about the sum of two probabilities, and the issues that motivate Cone Tightening or GMS go away. We verified that all testing methods successfully recover this and achieve excellent size control, including if tuning parameters are set to $0$.} This allows us to alternatively test rationalizability through a moment inequalities test that ensures uniform validity through GMS.\footnote{The implementation uses a ``Modified Method of Moments" criterion function, i.e. $S_1$ in the terminology of \citeasnoun{andrews-soares}, and the hard thresholding GMS function, i.e. studentized intercepts above $-\kappa_N$ were set to $0$ and all others to $-\infty$. The tuning parameter is set to $\kappa_N=\sqrt{ln(N_j)}$.}

\begin{table}
\begin{tabular}{r|rrrrrrr}
 & $\pi_0$ & $\pi_1$ & $\pi_2$ & $\pi_3$ & $\pi_4$ & $\pi_5$ & $\pi_6$ \\
 \hline
 $x_{1|1}$ & .181 & .190 & .2 & .240 & .3 & .167 & .152 \\
 $x_{2|1}$ & .226 & .238 & .25 & .213 & .2 & .167 & .107 \\
 $x_{3|1}$ & .226 & .238 & .25 & .213 & .2 & .333 & .441 \\
 $x_{4|1}$ & .367 & .333 & .3 & .333 & .3 & .333 & .3 \\
 $x_{1|2}$ & .181 & .190 & .2 & .240 & .3 & .333 & .486 \\
 $x_{2|2}$ & .226 & .238 & .25 & .213 & .2 & .167 & .107 \\
 $x_{3|2}$ & .226 & .238 & .25 & .213 & .2 & .167 & .107 \\
 $x_{4|2}$ & .367 & .333 & .3 & .333 & .3 & .333 & .3 \\
 $x_{1|3}$ & .181 & .190 & .2 & .240 & .3 & .333 & .486 \\
 $x_{2|3}$ & .226 & .238 & .25 & .213 & .2 & .167 & .107 \\
 $x_{3|3}$ & .226 & .238 & .25 & .213 & .2 & .167 & .107 \\
 $x_{4|3}$ & .367 & .333 & .3 & .333 & .3 & .333 & .3 \\
 \multicolumn{8}{c}{} \\
\end{tabular}
\caption{The $\pi$-vectors used for the Monte Carlo simulations in Table \ref{table:MC-2}.}
\label{table:pi-vectors}
\end{table}

\begin{table}
\begin{tabular}{l|l|rrrrrrrrrrr}
Method & $N_j$ & $\pi_0$ & $\cdot$ & $\cdot$ & $\cdot$ & $\cdot$ & $\pi_1$ & $\cdot$ & $\cdot$ & $\cdot$ & $\cdot$ & $\pi_2$ \\
\hline
Cone & 100 & .002 & .004 & .008 & .008 & .018 & .024 & .060 & .110 & .178 & .238 & .334 \\
Tightening & 200 & 0 & 0 & .004 & .008 & .012 & .040 & .088 & .164 & .286 & .410 & .544 \\
& 500 & 0 & 0 & .004 & .004 & .026 & .066 & .166 & .310 & .500 & .690 & .856 \\
& 1000 & 0 & 0 & 0 & 0 & .010 & .058 & .206 & .466 & .764 & .924 & .984 \\
\hline
GMS & 100 & 0 & 0 & 0 & .002 & .004 & .004 & .010 & .024 & .032 & .066 & .100 \\
& 200 & 0 & 0 & 0 & 0 & .004 & .010 & .018 & .048 & .082 & .174 & .262 \\
& 500 & 0 & 0 & .002 & .002 & .026 & .050 & .150 & .266 & .442 & .640 & .814 \\
& 1000 & 0 & 0 & 0 & 0 & .008 & .050 & .194 & .460 & .758 & .924 & .976 \\
\end{tabular}
\bigskip

\begin{tabular}{l|l|rrrrrrrrrrr}
Method & $N_j$ & $\pi_0$ & $\cdot$ & $\cdot$ & $\cdot$ & $\cdot$ & $\pi_3$ & $\cdot$ & $\cdot$ & $\cdot$ & $\cdot$ & $\pi_4$ \\
\hline
Cone & 100 & .002 & .002 & .008 & .006 & .012 & .016 & .036 & .070 & .098 & .148 & .200 \\
Tightening & 200 & 0 & 0 & .004 & .008 & .010 & .036 & .068 & .112 & .194 & .296 & .404  \\
& 500 &  0 & 0 & .004 & .004 & .026 & .064 & .156 & .296 & .456 & .664 & .786 \\
& 1000 & 0 & 0 & 0 & 0 & .010 & .058 & .200 & .460 & .756 & .916 & .974 \\
\hline
GMS & 100 & 0 & 0 & 0 & .002 & .004 & .004 & .006 & .018 & .028 & .066 & .120 \\
& 200 & 0 & 0 & 0 & 0 & .004 & .004 & .012 & .024 & .056 & .110 & .216 \\
& 500 & 0 & 0 & .002 & 0 & .004 & .012 & .034 & .084 & .196 & .348 & .544 \\
& 1000 & 0 & 0 & 0 & 0 & .002 & .028 & .080 & .180 & .406 & .694 & .892 \\
\end{tabular}
\bigskip

\begin{tabular}{l|l|rrrrrrrrrrr}
Method & $N_j$ & $\pi_0$ & $\cdot$ & $\cdot$ & $\cdot$ & $\cdot$ & $\pi_5$ & $\cdot$ & $\cdot$ & $\cdot$ & $\cdot$ & $\pi_6$ \\
\hline
Cone & 100 & .002 & .002 & .006 & .004 & .020 & .052 & .126 & .326 & .548 & .766 & .934 \\
Tightening & 200 & 0 & 0 & .002 & .004 & .008 & .044 & .202 & .490 & .836 & .962 & .996 \\
& 500 & 0 & 0 & .002 & .004 & .012 & .072 & .374 & .880 & .992 & 1 & 1  \\
& 1000 & 0 & 0 & 0 & 0 & .006 & .052 & .606 & .992 & 1 & 1 & 1 \\
\hline
GMS & 100 & 0 & 0 & 0 & .004 & .008 & .022 & .072 & .170 & .388 & .604 & .808 \\
& 200 & 0 & 0 & 0 & 0 & .006 & .020 & .112 & .300 & .640 & .880 & .974 \\
& 500 & 0 & 0 & 0 & .002 & .004 & .072 & .292 & .716 & .952 & 1 & 1 \\
& 1000 & 0 & 0 & 0 & 0 & 0 & .042 & .498 & .956 & 1 & 1 & 1 \\
\end{tabular}
\bigskip

\caption{Monte Carlo results. See Table \ref{table:pi-vectors} for definition of $\pi$-vectors. Recall that $\{\pi_1,\pi_3,\pi_5\}$ are on the boundary of $\mathcal{C}$ and $\pi_0$ is interior to it. All entries computed from $500$ simulations and $499$ replications per bootstrap.}
\label{table:MC-2}
\end{table}

Data were generated from a total of 31 d.g.p.'s described below and for sample sizes of $N_j \in \{100,200,500,1000\}$; recall that these are per budget, i.e. each simulated data set is based on 3 such samples. The d.g.p.'s are parameterized by the $\pi$-vectors reported in Table \ref{table:pi-vectors}. They are related as follows: $\pi_0$ is in the interior of $\mathcal{C}$; $\pi_2$, $\pi_4$, and $\pi_6$ are outside it; and $\pi_2$, $\pi_3$, and $\pi_5$ are on its boundary. Furthermore, $\pi_1=(\pi_0+\pi_2)/2$, $\pi_3=(\pi_0+\pi_4)/2$, and $\pi_5=(\pi_0+\pi_6)/2$. Thus, the line segment connecting $\pi_0$ and $\pi_2$ intersects the boundary of $\mathcal{C}$ precisely at $\pi_1$ and similarly for the next two pairs of vectors. We compute ``power curves" along those $3$ line segments at $11$ equally spaced points, i.e. changing mixture weights in increments of $.1$. This is replicated $500$ times at a bootstrap size of $R=499$. Nominal size of the test is $\alpha=.05$ throughout. Ideally, it should be exactly attained at the vectors $\{\pi_1,\pi_3,\pi_5\}$.

Results are displayed in Table \ref{table:MC-2}. Noting that the vectors are not too different, we would argue that the simulations indicate reasonable power. Adjustments that ensure uniform validity of tests do tend to cause conservatism for both GMS and cone tightening, but size control markedly improves with sample size.\footnote{We attribute some very slight nonmonotonicities in the ``power curves" to simulation noise.} While Cone Tightening appears less conservative than GMS in these simulations, we caution that the tuning parameters and the distance metrics underlying the test statistics are not directly comparable.

The differential performance across the three families of d.g.p.'s is expected because the d.g.p.'s were designed to pose different challenges. For both $\pi_1$ and $\pi_3$, one constraint is binding and three more are close enough to binding that, at the relevant sample sizes, they cannot be ignored. This is more the case for $\pi_3$ compared to $\pi_1$. It means that GMS or Cone Tightening will be necessary, but also that they are expected to be conservative. The vector $\pi_5$ has three constraints binding, with two more somewhat close. This is a worst case for naive (not using Cone Tightening or GMS) inference, which will rarely pick up all binding constraints. Indeed, we verified that inference with $\tau_N=0$ or $\kappa_N=0$ leads to overrejection. Finally, $\pi_2$ and $\pi_4$ fulfill the necessary conditions identified by \citeasnoun{Kawaguchi}, so that his test will have no asymptotic power at a parameter value in the first two panels of Table \ref{table:MC-2}.

\section{Empirical Application}

\label{sec:empirics}

\label{sec:Empirical Application}
We apply our methods to data from the U.K. Family Expenditure Survey, the same data used by BBC. Our testing of a RUM can, therefore, be compared with their revealed preference analysis of a representative consumer. To facilitate this comparison, we use the same selection from these data, namely the time periods from 1975 through 1999 and households with a car and at least one child. The number of data points used varies from 715 (in 1997) to 1509 (in 1975), for a total of 26341. For each year, we extract the budget corresponding to that year's median expenditure and, following Section \ref{sec:extending}, estimate the distribution of demand on that budget with polynomials of order $3$. Like BBC, we assume that all consumers in one year face the same prices, and we use the same price data. While budgets have a tendency to move outward over time, there is substantial overlap of budgets at median expenditure. To account for endogenous expenditure, we again follow Section \ref{sec:extending}, using total household income as instrument. This is also the same instrument used in BBC (2008).

We present results for blocks of eight consecutive periods and the same three composite goods (food, nondurable consumption goods, and services) considered in BBC.\footnote{As a reminder, Figure \ref{figure 2} illustrates the application. The budget is the 1993 one as embedded in the 1986-1993 block of periods, i.e. the figure corresponds to a row of Table \ref{table:8 periods}.} For all blocks of seven consecutive years, we analyze the same basket but also increase the dimensionality of commodity space to 4 or even 5. This is done by first splitting nondurables into clothing and other nondurables and then further into clothing, alcoholic beverages, and other nondurables. Thus, the separability assumptions that we (and others) implicitly invoke are successively relaxed. We are able to go further than much of the existing literature in this regard because, while computational expense increases with $K$, our approach is not subject to a statistical curse of dimensionality.\footnote{Tables \ref{table:7 periods} and \ref{table:8 periods} were computed in a few days on Cornell's ECCO cluster (32 nodes). An individual cell of a table can be computed in reasonable time on any desktop computer. Computation of a matrix $A$ took up to one hour and computation of one $\mathcal{J}_N$ about five seconds on a laptop.}

\begin{table}
\begin{tabular}{l|rrrr|rrrr|rrrr}
\multicolumn{1}{c}{} & \multicolumn{4}{c}{\textbf{3 goods}} & \multicolumn{4}{c}{\textbf{4 goods}}
& \multicolumn{4}{c}{\textbf{5 goods}} \\ 
& $I$ & $H$ & $\mathcal{J}_N$ & $p$ & $I$ & $H$ & $\mathcal{J}_N$ & $p$ & $I$ & $H$ & $\mathcal{J}_N$ & $p$ \\ 
\hline
75-81 & $36$ & $6409$ & $3.67$ & $.38$ & $52$ & $39957$ & $5.43$ & $.29$ & $55$ & $53816$ & $4.75$ & $.24$ \\ 
76-82 & $39$ & $4209$ & $11.6$ & $.14$ & $65$ & $82507$ & $5.75$ & $.39$ & $65$ & $82507$ & $5.34$ & $.31$ \\ 
77-83 & $41$ & $7137$ & $9.81$ & $.17$ & $65$ & $100728$ & $6.07$ & $.39$ & $68$ & $133746$ & $4.66$ & $.38$ \\ 
78-84 & $32$ & $3358$ & $7.38$ & $.24$ & $62$ & $85888$ & $2.14$ & $.70$ & $67$ & $116348$ & $1.45$ & $.71$ \\ 
79-85 & $35$ & $5628$ & $.114$ & $.96$ & $71$ & $202686$ & $.326$ & $.92$ & $79$ & $313440$ & $.219$ & $.94$ \\ 
80-86 & $38$ & $7104$ & $.0055$ & $.998$ & $58$ & $68738$ & $1.70$ & $.81$ & $66$ & $123462$ & $7.91$ & $.21$ \\ 
81-87 & $26$ & $713$ & $.0007$ & $.998$ & $42$ & $9621$ & $.640$ & $.89$ & $52$ & $28089$ & $6.33$ & $.27$ \\ 
82-88 & $15$ & $42$ & $0$ & $1$ & $21$ & $177$ & $.298$ & $.60$ & $31$ & $1283$ & $9.38$ & $.14$ \\ 
83-89 & $13$ & $14$ & $0$ & $1$ & $15$ & $31$ & $.263$ & $.49$ & $15$ & $31$ & $9.72$ & $.13$ \\ 
84-90 & $15$ & $42$ & $0$ & $1$ & $15$ & $42$ & $.251$ & $.74$ & $15$ & $42$ & $10.25$ & $.24$ \\ 
85-91 & $15$ & $63$ & $.062$ & $.77$ & $19$ & $195$ & $3.59$ & $.45$ & $21$ & $331$ & $3.59$ & $.44$ \\ 
86-92 & $24$ & $413$ & $1.92$ & $.71$ & $33$ & $1859$ & $7.27$ & $.35$ & $35$ & $3739$ & $9.46$ & $.28$ \\ 
87-93 & $45$ & $17880$ & $1.33$ & $.74$ & $57$ & $52316$ & $6.60$ & $.44$ & $70$ & $153388$ & $6.32$ & $.38$ \\ 
88-94 & $39$ & $4153$ & $1.44$ & $.70$ & $67$ & $136823$ & $6.95$ & $.38$ & $77$ & $313289$ & $6.91$ & $.38$ \\ 
89-95 & $26$ & $840$ & $.042$ & $.97$ & $69$ & $134323$ & $4.89$ & $.35$ & $78$ & $336467$ & $5.84$ & $.31$ \\ 
90-96 & $19$ & $120$ & $.040$ & $.95$ & $56$ & $52036$ & $4.42$ & $.19$ & $76$ & $272233$ & $3.55$ & $.25$ \\ 
91-97 & $17$ & $84$ & $.039$ & $.93$ & $40$ & $7379$ & $3.32$ & $.26$ & $50$ & $19000$ & $3.27$ & $.24$ \\ 
92-98 & $13$ & $21$ & $.041$ & $.97$ & $26$ & $897$ & $.060$ & $.93$ & $26$ & $897$ & $.011$ & $.99$ \\ 
93-99 & $9$ & $3$ & $.037$ & $.66$ & $15$ & $63$ & $0$ & $1$ & $15$ & $63$ & $0$ & $1$ \\
\multicolumn{13}{c}{} \\
\end{tabular}
\caption{Empirical results with 7 periods. $I=$ number of patches, $H=$ number of rationalizable discrete demand vectors, $\mathcal{J}_N=$ test statistic, $p=$ p-value.}
\label{table:7 periods}
\end{table}

\begin{table}
\begin{tabular}{l|rrrr}
\multicolumn{1}{c}{} & \multicolumn{4}{c}{\textbf{3 goods}} \\ 
& $I$ & $H$ & $\mathcal{J}_N$ & $p$ \\ 
\hline
75-82 & $51$ & $71853$ & $11.4$ & $.17$ \\ 
76-83 & $64$ & $114550$ & $9.66$ & $.24$ \\ 
77-84 & $52$ & $57666$ & $9.85$ & $.20$ \\ 
78-85 & $49$ & $76746$ & $7.52$ & $.26$ \\ 
79-86 & $55$ & $112449$ & $.114$ & $.998$ \\ 
80-87 & $41$ & $13206$ & $3.58$ & $.58$ \\ 
81-88 & $27$ & $713$ & $0$ & $1$ \\ 
82-89 & $16$ & $42$ & $0$ & $1$ \\ 
83-90 & $16$ & $42$ & $0$ & $1$ \\ 
84-91 & $20$ & $294$ & $.072$ & $.89$ \\ 
85-92 & $27$ & $1239$ & $2.24$ & $.68$ \\ 
86-93 & $46$ & $17880$ & $1.54$ & $.75$ \\ 
87-94 & $48$ & $39913$ & $1.55$ & $.75$ \\ 
88-95 & $42$ & $12459$ & $1.68$ & $.70$ \\ 
89-96 & $27$ & $840$ & $.047$ & $.97$ \\ 
90-96 & $24$ & $441$ & $.389$ & $.83$ \\ 
91-98 & $22$ & $258$ & $1.27$ & $.52$ \\ 
92-99 & $14$ & $21$ & $.047$ & $.96$ \\
\multicolumn{5}{c}{} \\
\end{tabular}
\caption{Empirical results with 8 periods. $I=$ number of patches, $H=$ number of rationalizable discrete demand vectors, $\mathcal{J}_N=$ test statistic, $p=$ p-value.}
\label{table:8 periods}
\end{table}

Regarding the test's statistical power, increasing the dimensionality of commodity space can in principle cut both ways. The number of rationality constraints increases, and this helps if some of the new constraints are violated but adds noise otherwise. Also, the maintained assumptions become weaker: In principle, a rejection of stochastic rationalizability at 3 but not 4 goods might just indicate a failure of separability. 

Tables \ref{table:7 periods} and \ref{table:8 periods} summarize our empirical findings. They display test statistics, p-values, and the numbers $I$ of patches and $H$ of rationalizable demand vectors; thus, matrices $A$ are of size $(I \times H)$.  All entries that show $\mathcal{J}_N=0$ and a corresponding p-value of $1$ were verified to be true zeros, i.e. $\hat{\pi}_{EC}$ is rationalizable. All in all, it turns out that estimated choice probabilities are typically not stochastically rationalizable, but also that this rejection is not statistically significant.\footnote{In additional analyses not presented here, we replicated these tables using polynomials of degree 2, as well as setting $\tau_N=0$. The qualitative finding of many positive but insignificant test statistics remains. In isolation, this finding may raise questions about the test's power. However, the test exhibits reasonable power in our Monte Carlo exercise and also rejects rationalizability in an empirical application elsewhere \cite{Hubner}. \newline We also checked whether small but positive test statistics are caused by adding-up constraints, i.e. by the fact that all components of $\pi$ that correspond to one budget must jointly be on some unit simplex. The estimator $\hat{\pi}$ can slightly violate this. Adding-up failures occur but are at least one order of magnitude smaller than the distance from a typical $\hat{\pi}$ to the corresponding projection $\hat{\eta}$.}

We identified a mechanism that may explain this phenomenon. Consider the 84-91 entry in Table \ref{table:8 periods}, where $\mathcal{J}_N$ is especially low. It turns out that one patch on budget $\mathcal{B}_5$ is below $\mathcal{B}_8$ and two patches on $\mathcal{B}_8$ are below $\mathcal{B}_5$. By the reasoning of Example \ref{example 1}, probabilities of these patches must add to less than $1$. The estimated sum equals $1.006$, leading to a tiny and statistically insignificant violation. This phenomenon occurs frequently and seems to cause the many positive but insignificant values of $\mathcal{J}_N$. The frequency of its occurrence, in turn, has a simple cause that may also appear in other data: If two budgets are slight rotations of each other and demand distributions change continuously in response, then population probabilities of patches like the above will sum to \textit{just} less than $1$. If these probabilities are estimated independently across budgets, the estimates will frequently add to slightly more than $1$. With $7$ or $8$ mutually intersecting budgets, there are many opportunities for such reversals, and positive but insignificant test statistics may become ubiquitous.

The phenomenon of estimated choice frequencies typically not being rationalizable means that there is need for a statistical testing theory and also a theory of rationality constrained estimation. The former is this paper's main contribution. We leave the latter for future research.

\section{Conclusion}

\label{sec:conclusion}

This paper presented asymptotic theory and computational tools for nonparametric testing of Random Utility Models. Again, the null to be tested was that data was generated by a RUM, interpreted as describing a heterogeneous population, where the only restrictions imposed on individuals' behavior were \textquotedblleft more is better\textquotedblright\ and SARP. In particular, we allowed for unrestricted, unobserved heterogeneity and stopped far short of assumptions that would recover invertibility of demand. We showed that testing the model is nonetheless possible. The method is easily adapted to choice problems that are discrete to begin with, and one can easily impose more, or fewer, restrictions at the individual level.

Possibilities for extensions and refinements abound, and some of these have already been explored. We close by mentioning further salient issues.\smallskip

\noindent \textbf{(1)} We provide algorithms (and code) that work for reasonably sized problem, but it would be extremely useful to make further improvements in this dimension.\smallskip

\noindent \textbf{(2)} The extension to infinitely many budgets is of obvious interest. Theoretically, it can be handled by considering an appropriate discretization argument \cite{mcfadden-2005}. For the proposed projection-based econometric methodology, such an extension requires evaluating choice probabilities locally over points in the space of $p$ via nonparametric smoothing, then use the choice probability estimators in the calculation of the $\mathcal{J}_N$-statistic. The asymptotic theory then needs to be modified. Another approach that can mitigate the computational constraint is to consider a partition of the space of $p$ such that $\mathbf{R}_{+}^{K}=\mathcal{P}_{1}\cup \mathcal{P}_{2}\cdots \cup \mathcal{P}_{M}$. Suppose we calculate the $\mathcal{J}_N$-statistic for each of these partitions. Given the resulting $M$-statistics, say $\mathcal{J}_N^{1},\cdots ,\mathcal{J}_N^{M}$, we can consider $\mathcal{J}_N^{\mathrm{max}}:=\max_{1\leq m\leq M}\mathcal{J}_N^{m}$ or a weighted average of them. These extensions and their formal statistical analysis are of practical interest.\smallskip

\noindent \textbf{(3)} While we allow for endogenous expenditure and therefore for the distribution of $u$ to vary with observed expenditure, we do assume that samples for all budgets are drawn from the same underlying population. This assumption can obviously not be dropped completely. However, it will frequently be of interest to impose it only conditionally on observable covariates, which must then be controlled for. This may be especially relevant for cases where different budgets correspond to independent markets, but also to adjust for slow demographic change as does, strictly speaking, occur in our data. It requires incorporating nonparametric smoothing in estimating choice probabilities as in Section \ref{sec:extending}, then averaging the corresponding $\mathcal{J}_N$-statistics over the covariates. This extension will be pursued.\smallskip

\noindent \textbf{(4)} Natural next steps after rationality testing are extrapolation to (bounds on) counterfactual demand distributions and welfare analysis, i.e. along the lines of BBC (2008) or, closer to our own setting, \citeasnoun{Adams16} and \citeasnoun{DKQS16}. This extension is being pursued. Indeed, the tools from Section \ref{sec:population} have already been used for choice extrapolation (using algorithms from an earlier version of this paper) in \citeasnoun{Manski14}.  \smallskip

\noindent \textbf{(5)} The econometric techniques proposed here can be potentially useful in much broader contexts. Indeed, they have already been used to nonparametrically test game theoretic model with strategic complementarities \cite{LQS15,LQS18}, a novel model of ``price preference'' \cite{DKQS16}, and the collective household model \cite{Hubner}. Even more generally, existing proposals for testing in moment inequality models \cite{andrews-guggenberger-2009,andrews-soares,BCS15,romano-shaikh} work with explicit inequality constraints, i.e. (in the linear case) $\mathcal{H}$-representations. In settings in which theoretical restrictions inform a $\mathcal{V}$-representation of a cone or, more generally, a polyhedron, the $\mathcal{H}$-representation will typically not be available in practice. We expect that our method can be used in many such cases.

\nocite{BBC03}
\nocite{BBC07}
\nocite{BBC08}

\newpage

\bibliographystyle{econometrica}
\bibliography{mybib}

\clearpage
\pagenumbering{arabic}
\renewcommand*{\thepage}{S\arabic{page}}
\renewcommand{\theequation}{S.\arabic{equation}}
\renewcommand{\thecondition}{S.\arabic{condition}}

\begin{center}
	{\large \bf  SUPPLEMENTAL MATERIALS}
\end{center}

\setcounter{section}{0}

\section*{Appendix A: Proofs and Further Details of Inferential Procedures}
\label{appendix A}

\begin{proof}[\textup{\textbf{Proof of Theorem~\protect\ref{prop:cone}}}]
The proof uses nonstochastic demand systems, which can be identified with vectors $(d_1,\dots,d_J) \in \mathcal{B}_1 \times \dots \times \mathcal{B}_J$. Such a system is rationalizable if $d_j \in \arg \max_{y \in \mathcal{B}_j}{u(y)}, j=1,\dots,J$ for some utility function $u$.

Rationalizability of nonstochastic demand systems is well understood. In particular, and irrespective of whether we define rationalizability by GARP or SARP, it is decidable from knowing the preferences directly revealed by choices, hence from knowing patches containing $(d_1,\dots,d_J)$. It follows that for all nonstochastic demand systems that select from the same patches, either all or none are rationalizable.

Fix $(P_1,...,P_J)$. Let the set $\mathcal{Y}^*$ collect one ``representative" element (e.g., the geometric center point) of each patch. Let $(P^*_1,...,P^*_J)$ be the unique stochastic demand system concentrated on $\mathcal{Y}^*$ and having the same vector representation as $(P_1,...,P_J)$. The previous paragraph established that demand systems can be arbitrarily perturbed within patches, so $(P_1,...,P_J)$ is rationalizable iff $(P^*_1,...,P^*_J)$ is. It follows that rationalizability of $(P_1,...,P_J)$ can be decided from its vector representation $\pi$, and that it suffices to analyze stochastic demand systems supported on $\mathcal{Y}^*$. 

Now, any stochastic demand system is rationalizable iff it is a mixture of rationalizable nonstochastic ones. Since $\mathcal{Y}^*$ is finite, there are finitely many nonstochastic demand systems supported on it; of these, a subset will be rationalizable. Noting that these demand systems are characterized by binary vector representations corresponding to columns of $A$, the statement of the Theorem is immediate for the restricted class of stochastic demand systems supported on $\mathcal{Y}^*$.
\end{proof}

\begin{proof}[\textup{\textbf{Proof of Theorem~\protect\ref{prop1}}}]
We begin with some preliminary observations. Throughout this proof, $c(\mathcal{B}_{i}) $ denotes the object actually chosen from budget $\mathcal{B}_{i}$.

(i) If there is a choice cycle of any finite length, then there is a cycle of length 2 or 3 (where a cycle of length 2 is a WARP violation). To see this, assume there exists a length $N$ choice cycle $c(\mathcal{B}_{i})\succ c(\mathcal{B}_{j})\succ c(\mathcal{B}_{k})\succ ...\succ c(\mathcal{B}_{i})$. If $c(\mathcal{B}_{k})\succ c(\mathcal{B}_{i})$, then a length 3 cycle has been discovered. Else, there exists a length $N-1$ choice cycle $c(\mathcal{B}_{i})\succ c(\mathcal{B}_{k})\succ ...\succ c(\mathcal{B}_{i})$. The argument can be iterated until $N=4$.

(ii) Call a length 3 choice cycle \textit{irreducible} if it does not contain a length 2 cycle. Then a choice pattern is rationalizable iff it contains no length 2 cycles and also no irreducible length 3 cycles. (In particular, one can ignore reducible length 3 cycles.) This follows trivially from (i).

(iii) Let $J=3$ and $M=1$, i.e. assume there are three budgets but two of them fail to intersect. Then any length 3 cycle is reducible. To see this, assume w.l.o.g. that $\mathcal{B}_{1}$ is below $\mathcal{B}_{3}$, thus $c(\mathcal{B}_{3})\succ c(\mathcal{B}_{1})$ by monotonicity. If there is a choice cycle, we must have $c(\mathcal{B}_{1})\succ c(\mathcal{B}_{2})\succ c(\mathcal{B}_{3})$. $c(\mathcal{B}_{1})\succ c(\mathcal{B}_{2})$ implies that $c(\mathcal{B}_{2})$ is below $\mathcal{B}_{1}$, thus it is below $\mathcal{B}_{3}$. $c(\mathcal{B}_{2})\succ c(\mathcal{B}_{3})$ implies that $c(\mathcal{B}_{3})$ is below $\mathcal{B}_{2}$. Thus, choice from $(\mathcal{B}_{2},\mathcal{B}_{3})$ violates WARP.

We are now ready to prove the main result. The nontrivial direction is \textquotedblleft only if,\textquotedblright\ thus it suffices to show the following: If choice from $(\mathcal{B}_{1},...,\mathcal{B}_{J-1})$ is rationalizable but choice from $(\mathcal{B}_{1},...,\mathcal{B}_{J})$ is not, then choice from $(\mathcal{B}_{M+1},...,\mathcal{B}_{J})$ cannot be rationalizable. By observation (ii), if $(\mathcal{B}_{1},...,\mathcal{B}_{J})$ is not rationalizable, it contains either a 2-cycle or an irreducible 3-cycle. Because choice from all triplets within $(\mathcal{B}_{1},...,\mathcal{B}_{J-1})$ is rationalizable by assumption, it is either the case that some $(\mathcal{B}_{i},\mathcal{B}_{J})$ constitutes a 2-cycle or that some triplet $(\mathcal{B}_{i},\mathcal{B}_{k},\mathcal{B}_{J})$, where $i<k$ w.l.o.g., reveals an irreducible choice cycle. In the former case, $\mathcal{B}_{i}$ must intersect $\mathcal{B}_{J}$, hence $i>M$, hence the conclusion. In the latter case, if $k\leq M$, the choice cycle must be a 2-cycle in $(\mathcal{B}_{i},\mathcal{B}_{k})$, contradicting rationalizability of $(\mathcal{B}_{1},...,\mathcal{B}_{J-1})$. If $i\leq M$, the choice cycle is reducible by (iii). Thus, $i>M$, hence the conclusion.
\end{proof}

\begin{proof}[\textup{\textbf{Proof of Lemma~\protect\ref{lem:c_tau}}}]
	Letting $\nu_\tau = \nu - (\tau/H) \mathbf{1}_H$ in $\mathcal{C}_\tau = \{A\nu|\nu \geq (\tau/H){\bf 1}_H \}$
	we have 
	\begin{eqnarray*}
		\mathcal{C}_\tau & = &\{A[\nu_\tau + (\tau/H) \mathbf{1}_H] | \nu_\tau \geq 0\} \\
		& = & \mathcal{C} \oplus (\tau/H) A \mathbf{1}_H \\
		& = & \{t: t - (\tau/H) A \mathbf{1}_H \in \mathcal{C}\}
	\end{eqnarray*}
	where $\oplus$ signifies Minkowski sum. Define 
	\begin{equation*}
		\phi = - B A {\mathbf{1}}_H/H.
	\end{equation*}
	Using the $\mathcal{H}$-representation of $\mathcal{C}$, 
	\begin{eqnarray*}
		\mathcal{C}_\tau & = & \{t: B(t - (\tau/H) A \mathbf{1}_H) \leq 0\} \\
		& = & \{t: B t \leq -\tau \phi\}.
	\end{eqnarray*}
	Note that the above definition of $\phi$ implies $\phi \in \mathrm{col}(B)$.
	Also define 
	\begin{eqnarray*}
		\Phi &:=& -BA \\
		&=& -\left[ 
		\begin{array}{c}
			b_1^{\prime } \\ 
			\vdots \\ 
			b_m^{\prime }%
		\end{array}
		\right] \left [ a_1, \cdots, a_H \right ] \\
		&=& \{\Phi_{kh}\}
	\end{eqnarray*}
	where $\Phi_{kh} = b_k^{\prime }a_h, 1 \leq k \leq m, 1 \leq h \leq H$ and
	let $e_h$ be the $h$-th standard unit vector in $\mathbf{R}^H$. Since $e_h
	\geq 0$, the $\mathcal{V}$-representation of $\mathcal{C}$ implies that $A e_h \in \mathcal{C}$,
	and thus 
	\begin{equation*}
		BA e_h \leq 0
	\end{equation*}
	by its $\mathcal{H}$-representation. Therefore 
	\begin{equation}  \label{eq:phi}
		\Phi_{kh} = - e_k^{\prime }BA e_h \geq 0, \quad 1 \leq k \leq m, 1 \leq h
		\leq H.
	\end{equation}
	\noindent But if $k \leq \bar m$, it cannot be that 
	\begin{equation*}
		a_j \in \{x: b_k^{\prime }x = 0\} \quad \text{ for all } j 
	\end{equation*}
	whereas 
	$$
	b_k'a_h = 0  
	$$ 
	holds for $\bar m + 1 \leq k \leq m, 1 \leq h \leq H$.  
	Therefore if $k \leq \bar m$, $\Phi_{kh} = b_k^{\prime }a_h$ is nonzero at least
	for one $h, 1 \leq h \leq H$, whereas if $k > \bar m$, $\Phi_{kh} = 0$ for every $h$.    Since \eqref{eq:phi} implies that all of $%
	\{\Phi_{kh}\}_{h=1}^H$ are non-negative, we conclude that 
	\begin{equation*}
		\phi_k = \frac 1 H \sum_{h=1}^H \Phi_{kh} > 0
	\end{equation*}
	for every $k \leq \bar m$ and $\phi_k = 0$ for every $k > \bar m$. We now have 
	\begin{equation*}
		\mathcal{C}_\tau = \{t: B t \leq - \tau \phi\}
	\end{equation*}
	where $\phi$ satisfies the stated properties {\it (i)} and {\it (ii)}.  
\end{proof}

Before we present the proof of Theorem \ref{thm1}, it is necessary to specify a class of distributions, to which we impose a mild condition that guarantees stable behavior of the statistic $\mathcal{J}_N$. To this end, we further specify the nature of each row of $B$. Recall that w.l.o.g. the first $\bar m$ rows of $B$ correspond to inequality constraints, whereas the rest of the rows represent equalities. Note that the $\bar m$ inequalities include nonnegativity constraints $\pi_{i|j} \geq 0, 1 \leq i \leq I_j, 1 \leq j \leq J$, represented by the row of $B$ consisting of a negative constant for the corresponding element and zeros otherwise.  Likewise, the identities that $\sum_{i=1}^{I_j}\pi_{i|j}$ is constant across $1 \leq j \leq J$ are included in the set of equality constraints.\footnote{If we impose the (redundant) restriction ${\bf 1}_H'\nu = 1$ in the definition of $\mathcal{C}$, then the corresponding equality restrictions would be $\sum_{i=1}^{I_j}\pi_{i|j} = 1$ for every $j$.}  We show in the proof that the presence of these ``definitional'' equalities/inequalities, which always hold by construction of $\hat \pi$, do not affect the asymptotic theory even when they are (close to) be binding.  Define $\mathcal K = \{1,...,m\}$, and let ${\mathcal K}^D$ be the set of indices for the rows of $B$ corresponding to the above nonnegativity constraints and the constant-sum constraints.  Let ${\mathcal K}^R = \mathcal K \setminus {\mathcal K}^D$, so that $b_k'\pi \leq 0$ represents an economic restriction if $k \in {\mathcal K}^R$.\footnote{In \eqref{eq:b-matrix}, $\mathcal{K}^R$ contains only the last row of the matrix.}   Recalling the choice vectors  $(d_{j|1},...,d_{j|N_j})$ are IID-distributed within each time period $j, 1 \leq j \leq J$, let $d_j$ denote the choice vector of a consumer facing budget $j$ (therefore w.l.o.g we can let  $d_j = d_{j|1}$).     Define $d = [d_1',...,d_J']'$,  a random $I$-vector of binary variables.   Note ${\mathrm{E}[d] = \pi}$.   Let 
\begin{eqnarray*}
	g &=& Bd
	\\
	&=& [g_1,...,g_m]'.
\end{eqnarray*}
With these definitions, consider the following requirement:
\begin{condition}
	\label{condition 1} For each  $k \in {\mathcal K}^R$,  {\rm{var}}$(g_k)  > 0$ and $\mathrm{E}[|g_k/\sqrt{\mathrm{var}(g_k)}|^{2+c_1}] < c_2$ hold, where $c_1$ and $c_2$ are positive constants. 
\end{condition}
\noindent 
This  type of condition is standard in the literature; see, for example, \citeasnoun{andrews-soares}. 

\begin{proof}[\textup{\textbf{Proof of Theorem~\protect\ref{thm1}}}]
	
	By applying the Minkowski-Weyl theorem and Lemma \ref{lem:c_tau} to $\mathcal{J}_N$
	and $\tilde{J}_{N}(\tau _{N})$, we see that our procedure is equivalent to
	comparing 
	\begin{equation*}
		\mathcal{J}_N=\min_{t\in \mathbf{R}^{I}:Bt\leq 0}N[\hat{\pi}-t]^{\prime }\Omega
		\lbrack \hat{\pi}-t]
	\end{equation*}%
	to the $1-\alpha $ quantile of the distribution of 
	\begin{equation*}
		\tilde{\mathcal J}_{N}=\min_{t\in \mathbf{R}^{I}:Bt\leq -\tau _{N}\phi }N[%
		\tilde{\eta}_{\tau _{N}}-t]^{\prime }\Omega \lbrack \tilde{\eta}_{\tau
			_{N}}-t]
	\end{equation*}%
	with $\phi = [\bar \phi',(0,...,0)']'$,   $\bar \phi \in \mathbf{R}_{++}^{\bar m}$, where 
	\begin{equation*}
		\tilde{\eta}_{\tau _{N}}=\hat{\eta}_{\tau _{N}}+\frac{1}{\sqrt{N}}N(0,\hat{S}%
		),
	\end{equation*}%
	\begin{equation*}
		\hat{\eta}_{\tau _{N}}=\operatornamewithlimits{argmin}_{t\in \mathbf{R}%
			^{I}:Bt\leq -\tau _{N}\phi }N[\hat{\pi}-t]^{\prime }\Omega \lbrack \hat{\pi}%
		-t].
	\end{equation*}%
	Suppose $B$ has $m$ rows and rank$(B)=\ell $. Define an $\ell \times m$
	matrix $K$ such that $KB$ is a matrix whose rows consist of a basis of the
	row space row$(B)$. Also let $M$ be an $(I-\ell )\times I$ matrix whose rows
	form an orthonormal basis of ker$B=$ ker$(KB)$, and define $P=\binom{KB}{M}$%
	. Finally, let $\hat{g}=B\hat{\pi}$ and $\hat{h}=M\hat{\pi}$. Then 
	\begin{eqnarray*}
		\mathcal{J}_N &=&\min_{Bt\leq 0}N\left[ \binom{KB}{M}(\hat{\pi}-t)\right] ^{\prime }{%
			P^{-1}}^{\prime }\Omega P^{-1}\left[ \binom{KB}{M}(\hat{\pi}-t)\right] \\
		&=&\min_{Bt\leq 0}N\binom{K[\hat{g}-Bt]}{\hat{h}-Mt}^{\prime }{P^{-1}}%
		^{\prime }\Omega P^{-1}\binom{K[\hat{g}-Bt]}{\hat{h}-Mt}.
	\end{eqnarray*}%
	Let 
	\begin{equation*}
		{\mathcal{U}}_{1}=\left\{ \binom{K\gamma }{h}:\gamma =Bt,h=Mt,B^\leq t\leq 0,B^= t = 0, t\in 
		\mathbf{R}^{I}\right\}
	\end{equation*}%
	then writing $\alpha =KBt$ and $h=Mt$, 
	\begin{equation*}
		\mathcal{J}_N=\min_{\binom{\alpha }{h}\in \mathcal{U}_{1}}N\binom{K\hat{g}-\alpha }{%
			\hat{h}-h}^{\prime }{P^{-1}}^{\prime }\Omega P^{-1}\binom{K\hat{g}-\alpha }{%
			\hat{h}-h}.
	\end{equation*}%
	Also define 
	\begin{equation*}
		\mathcal{U}_{2}=\left\{ \binom{K\gamma }{h}:\gamma = \binom{\gamma^\leq}{\gamma^=}, \gamma^\leq \in {\bf R}^{\bar m}_+, \gamma^= = 0, \gamma \in \text{%
			col}(B),h\in \mathbf{R}^{I-\ell }\right\}
	\end{equation*}%
	where col$(B)$ denotes the column space of $B$. Obviously ${\mathcal{U}}%
	_{1}\subset {\mathcal{U}}_{2}$. Moreover, ${\mathcal{U}}_{2}\subset {%
		\mathcal{U}}_{1}$ holds. To see this, let $\binom{K\gamma ^{\ast }}{h^{\ast }%
	}$ be an arbitrary element of ${\mathcal{U}}_{2}$. We can always find $%
	t^{\ast }\in \mathbf{R}^{I}$ such that $\gamma ^{\ast }=Bt^{\ast }$.
	Define 
	\begin{equation*}
		t^{\ast \ast }:=t^{\ast }+M^{\prime }h^{\ast }-M^{\prime }Mt^{\ast }
	\end{equation*}%
	then $Bt^{\ast \ast }=Bt^{\ast }=\gamma ^{\ast }$, therefore $B^\leq t^{**} \leq 0$ and $B^= t^{**} = 0$.  Also, $Mt^{\ast \ast
	}=Mt^{\ast }+MM^{\prime }h^{\ast }-MM^{\prime }Mt^{\ast }=h^{\ast }$,
	therefore $\binom{K\gamma ^{\ast }}{h^{\ast }}$ is an element of ${\mathcal{U%
		}}_{1}$ as well. Consequently, 
		\begin{equation*}
			{\mathcal{U}}_{1}={\mathcal{U}}_{2}.
		\end{equation*}%
		We now have 
		\begin{eqnarray*}
			\mathcal{J}_N &=&\min_{\binom{\alpha }{h}\in \mathcal{U}_{2}}N\binom{K\hat{g}-\alpha 
			}{\hat{h}-h}^{\prime }{P^{-1}}^{\prime }\Omega P^{-1}\binom{K\hat{g}-\alpha 
		}{\hat{h}-h} \\
		&=&N\min_{\binom{\alpha }{y}\in \mathcal{U}_{2}}\binom{K\hat{g}-\alpha }{y}%
		^{\prime }{P^{-1}}^{\prime }\Omega P^{-1}\binom{K\hat{g}-\alpha }{y}.
	\end{eqnarray*}%
	Define 
	\begin{equation*}
		T(x,y)=\binom{x}{y}^{\prime }{P^{-1}}^{\prime }\Omega P^{-1}\binom{x}{y}%
		,\quad x\in \mathbf{R}^{\ell },y\in \mathbf{R}^{I-\ell },
	\end{equation*}%
	and 
	\begin{equation*}
		t(x):=\min_{y\in \mathbf{R}^{I-\ell }}T(x,y),\quad s(g):=
		\min_{\gamma = [{\gamma^\leq}',{\gamma^=}']', \gamma^\leq \leq 0, \gamma^= = 0, \gamma \in \text{col}(B)}t(K[g-\gamma ]).
	\end{equation*}%
	It is easy to see that $t:\mathbf{R}^{\ell }\rightarrow \mathbf{R}_{+}$ is a
	positive definite quadratic form. We can write 
	\begin{eqnarray*}
		\mathcal{J}_N &=&N\min_{\gamma = [{\gamma^\leq}',{\gamma^=}']', \gamma^\leq \leq 0, \gamma^= = 0, \gamma \in \text{col}(B)}t(K[\hat{g}-\gamma ])
		\\
		&=&Ns(\hat{g})\quad \quad \quad \\
		&=&s(\sqrt{N}\hat{g}).
	\end{eqnarray*}%
	We now show that tightening can turn non-binding inequality constraints into
	binding ones but not vice versa. Note that, as will be seen below, this
	observation uses diagonality of $\Omega$ and the specific geometry of the
	cone $\mathcal{C}$. Let $\hat{\gamma}_{\tau _{N}}^{k}$, $\hat{g}^{k}$ and $\phi ^{k}$
	denote the $k$-th elements of $\hat{\gamma}_{\tau _{N}}=B\hat{\eta}_{\tau
		_{N}}$, $\hat{g}$ and $\phi $. Moreover, define $\gamma _{\tau }(g)=[\gamma
	^{1}(g),...,\gamma ^{m}(g)]^{\prime }=\operatornamewithlimits{argmin}_{\gamma = [{\gamma^\leq}',{\gamma^=}']', \gamma^\leq \leq -\tau\bar \phi, \gamma^= = 0, \gamma \in \text{col}(B)}
	t(K[g-\gamma ])$ for $%
	g\in {\mathrm{col}(B)}$, and let $\gamma_\tau^k(g)$ be its $k$-th element.  Then $\hat{\gamma}_{\tau _{N}}=\gamma _{\tau _{N}}(%
	\hat{g})$. Finally, define $\beta _{\tau }(g)=\gamma _{\tau }(g)+\tau \phi $
	for $\tau >0$ and let $\beta _{\tau }^{k}(g)$ denote its $k$-th element.  Note $\gamma_\tau^k(g) = \phi^k = \beta _{\tau }^{k}(g)=0$ for every $k > \bar m$ and $g$.   Now
	we show that for each $k \leq \bar m$ and for some $\delta >0$, 
	\begin{equation*}
		\beta _{\tau }^{k}(g)=0
	\end{equation*}%
	if $|g^{k}|\leq \tau \delta $ and $g^{j}\leq \tau \delta ,1\leq j\leq \bar m$. In
	what follows we first show this for the case with $\Omega =\mathbf{I}_{I}$,
	where $\mathbf{I}_{I}$ denotes the $I$-dimensional identity matrix, then
	generalize the result to the case where $\Omega $ can have arbitrary
	positive diagonal elements.
	
	For $\tau >0$ and $\delta >0$ define hyperplanes 
	\begin{equation*}
		H_{k}^{\tau }=\{x:b_{k}^{\prime }x=-\tau \phi ^{k}\},
	\end{equation*}
	\begin{equation*}
		H_{k}=\{x:b_{k}^{\prime }x=0\},
	\end{equation*}%
	half spaces%
	\begin{equation*}
		{H_{\angle }^{\tau }}_{k}(\delta )=\{x:b_{k}^{\prime }x\leq \tau \delta
		\},
	\end{equation*}%
	and also  
	\begin{equation*}
		S_{k}(\delta )=\{x\in \mathcal{C}:|b_{k}^{\prime }x|\leq \tau \delta \}
	\end{equation*}%
	for $1\leq k\leq m$. Define 
	$$
	L = \cap_{k=\bar m + 1}^m H_{k},
	$$
	a linear subspace of ${\bf R}^I$.  In what follows we show that for small
	enough $\delta >0$, every element $x^{\ast }\in \mathbf{R}^{I}$ such that 
	\begin{equation}
		x^{\ast }\in S_{1}(\delta )\cap \cdots \cap S_{q}(\delta )\cap {H_{\angle
			}^{\tau }}_{q+1}(\delta )\cap \cdots {H_{\angle }^{\tau }}_{m}(\delta )\text{
			for some }q\in \{1,...,\bar m\}  \label{x^*}
	\end{equation}%
	satisfies 
	\begin{equation}
		x^{\ast }|\mathcal{C}_{\tau }\in H_{1}^{\tau }\cap \cdots \cap H_{q}^{\tau }\cap L
		\label{projx^*}
	\end{equation}%
	where $x^{\ast }|\mathcal{C}_{\tau }$ denotes the orthogonal projection of $x^{\ast }$
	on $\mathcal{C}_{\tau }$. Let ${g^{\ast }}^{k}=b_{k}^{\prime }x^{\ast },k=1,...,m$.
	Note that an element $x^{\ast }$ fulfills \eqref{x^*} iff $|{g^{\ast }}%
	^{k}|\leq \tau \delta ,1\leq k\leq q$ and ${g^{\ast }}^{j}\leq \tau \delta
	,q+1\leq j\leq \bar m$. Likewise, \eqref{projx^*} holds iff $\beta _{k}^{\tau
	}(g^{\ast })=0,1\leq k\leq q$ (recall $\beta _{k}^{\tau
}(g^{\ast })=0$ always holds for $k > \bar m$). Thus in order to establish the desired
property of the function $\beta _{\tau }(\cdot )$, we show that \eqref{x^*}
implies \eqref{projx^*}. Suppose it does not hold; then without loss of
generality, for an element $x^{\ast }$ that satisfies \eqref{x^*} for an
arbitrary small $\delta >0$, we have 
\begin{equation}
	x^{\ast }|\mathcal{C}_{\tau }\in H_{1}^{\tau }\cap \cdots \cap H_{r}^{\tau }\cap L\quad 
	\text{ and }\quad x^{\ast }|\mathcal{C}_{\tau }\notin H_{j}^{\tau },r+1\leq j\leq q
	\label{xandc}
\end{equation}%
for some $1\leq r\leq q-1$. Define halfspaces 
\begin{equation*}
	{H_{\angle }^{\tau }}_{k}=\{x:b_{k}^{\prime }x\leq -\tau \phi ^{k}\},
\end{equation*}%
\begin{equation*}
	{H_{\angle }}_{k}=\{x:b_{k}^{\prime }x\leq 0\}
\end{equation*}%
for $1\leq k\leq m,\tau >0$ and also let 
\begin{equation*}
	F=H_{1}\cap \cdots \cap H_{r} \cap \mathcal{C},
\end{equation*}%
then for \eqref{xandc} to hold for some $x^{\ast }\in \mathbf{R}^{I}$
satisfying \eqref{x^*} for an arbitrary small $\delta >0$ we must have 
\begin{equation*}
	F|\left( H_{1}^{\tau }\cap \cdots \cap H_{r}^{\tau } \cap L\right) \subset \mathrm{%
		int}({H_{\angle }^{\tau }}_{r+1}\cap \cdots \cap {H_{\angle }^{\tau }}_{q})
\end{equation*}%
(Recall the notation $|$ signifies orthogonal projection. Also note that if
dim$(F) = 1$, then \eqref{xandc} does
not occur under \eqref{x^*}.) Therefore if we let 
\begin{equation*}
	\Delta (J)=\{x\in \mathbf{R}^{I}:\mathbf{1}_{I}^{\prime }x=J,x\geq 0\},
\end{equation*}%
i.e. the simplex with vertices $(J,0,\cdots ,0),\cdots ,(0,\cdots ,0,J)$, we
have 
\begin{equation}
	\left( F\cap \Delta (J)\right) |\left( H_{1}^{\tau }\cap \cdots \cap
	H_{r}^{\tau }\cap L\right) \subset \mathrm{int}({H_{\angle }^{\tau }}_{r+1}\cap
	\cdots \cap {H_{\angle }^{\tau }}_{q}).  \label{stickout}
\end{equation}%
Let $\{a_{1},...,a_{H}\}=\mathcal{A}$ denote the collection of the column
vectors of $A$. Then $\{\text{the vertices of }F\cap \Delta (J)\}\in 
\mathcal{A}$. Let $\bar{a},\bar{\bar{a}}\in F\cap \Delta (J)$. Let $%
B(\varepsilon ,x)$ denote the $\varepsilon $-(open) ball with center $x\in {%
	\mathbf{R}}^{I}$. By \eqref{stickout}, 
\begin{equation*}
	B\left( \varepsilon ,\left( \bar{a}|\cap _{j=1}^{r}H_{j}^{\tau }\cap L\right)
	\right) \subset \mathrm{int}({H_{\angle }^{\tau }}_{r+1}\cap \cdots \cap {%
		H_{\angle }^{\tau }}_{q})\cap {H_{\angle }}_{1}\cap \cdots \cap {H_{\angle }}%
	_{r}
\end{equation*}%
holds for small enough $\varepsilon >0$. Let $\bar{a}^{\tau }:=\bar{a}+\frac{\tau}{H}A {\bf 1}_H $%
, $\bar{\bar{a}}^{\tau }:=\bar{\bar{a}}+\frac{\tau}{H}A {\bf 1}_H$, then 
\begin{eqnarray*}
	\left( \left( \bar{a}|(\cap _{j=1}^{r}H_{j}^{\tau })\cap L\right) -\bar{a}\right)
	^{\prime }\left( \bar{\bar{a}}-\bar{a}\right) &=&\left( \left( \bar{a}|(\cap
	_{j=1}^{r}H_{j}^{\tau })\cap L\right) -\bar{a}\right) ^{\prime }\left( \bar{\bar{a}}%
	^{\tau }-\bar{a}^{\tau }\right) \\
	&=&0
\end{eqnarray*}%
since $\bar{a}^{\tau },\bar{\bar{a}}^{\tau }\in (\cap _{j=1}^{r}H_{j}^{\tau })\cap L$%
. We can then take $z\in B\left( \varepsilon ,\left( \bar{a}|(\cap
_{j=1}^{r}H_{j}^{\tau })\cap L\right) \right) $ such that $(z-\bar{a})^{\prime }(%
\bar{\bar{a}}-\bar{a})<0$. By construction $z\in \mathcal{C}$, which implies the
existence of a triplet $(a,\bar{a},\bar{\bar{a}})$ of distinct elements in $%
\mathcal{A}$ such that $(a-\bar{a})^{\prime }(\bar{\bar{a}}-\bar{a})<0$. In
what follows we show that this cannot happen, then the desired property of $%
\beta _{\tau }$ is established.

So let us now show that 
\begin{equation}  \label{acute}
	(a_1 - a_0)^{\prime }(a_2 - a_0) \geq 0 \text{ for every triplet }
	(a_0,a_1,a_2) \text{ of distinct elements in } \mathcal{A}.
\end{equation}
Noting that $a_{i}^{\prime }a_{j}$ just counts the number of budgets on
which $i$ and $j$ agree, define 
\begin{equation*}
	\phi(a_i,a_j)=J-a_{i}^{\prime }a_{j},
\end{equation*}
the number of disagreements. Importantly, note that $\phi(a_i,a_j)=
\phi(a_j,a_i)$ and that $\phi $ is a distance (it is the taxicab distance
between elements in $\mathcal{A}$, which are all 0-1 vectors). Now 
\begin{eqnarray*}
	&&(a_1 - a_0)^{\prime }(a_2 - a_0) \\
	&& =a_{1}^{\prime }a_{2}-a_{0}^{\prime }a_{2}-a_{1}^{\prime
	}a_{0}+a_{0}^{\prime }a_{0} \\
	&& =J-\phi(a_1,a_2)-(J-\phi(a_0,a_2))-(J-\phi(a_0,a_1))+J \\
	&& =\phi(a_0,a_2)+\phi(a_0,a_1)-\phi(a_1,a_2)\geq 0
\end{eqnarray*}%
by the triangle inequality.

Next we treat the case where $\Omega$ is not necessarily $\mathbf{I}_I$.
Write 
\begin{equation*}
	\Omega=\left[ 
	\begin{array}{cccc}
		\omega_1^2 & 0 & \ldots & 0 \\ 
		0 & \omega_2^2 & \ldots & 0 \\ 
		&  & \ddots &  \\ 
		0 & \ldots & 0 & \omega_I^2%
	\end{array}%
	\right].
\end{equation*}%
The statistic $\mathcal{J}_N$ in \eqref{J-def} can be rewritten, using the square-root
matrix $\Omega^{1/2}$, 
\begin{equation*}
	\mathcal{J}_N = \min_{\eta^* = \Omega^{1/2}\eta: \eta \in C}[\hat \pi^* -
	\eta^*]^{\prime }[\hat \pi^* - \eta^*]
\end{equation*}
or 
\begin{equation*}
	\mathcal{J}_N = \min_{\eta^* \in C^*}[\hat \pi^* - \eta^*]^{\prime }[\hat \pi^* -
	\eta^*]
\end{equation*}
where 
\begin{eqnarray*}
	\mathcal{C}^* &=& \{\Omega^{1/2}A\nu| \nu \geq 0 \} \\
	&=& \{A^*\nu|\nu \geq 0\}
\end{eqnarray*}
with 
\begin{equation*}
	A^* = [a^*_1,...,a^*_H], a^*_h = \Omega^{1/2}a_h, 1 \leq h \leq H.
\end{equation*}
Then we can follow our previous argument replacing $a$'s with $a^*$'s, and
using 
\begin{equation*}
	\Delta^*(J) = {\mathrm{conv}}([0,...,\omega_i,....0]^{\prime } \in \mathbf{R}%
	^I, i = 1,...,I).
\end{equation*}
instead of the simplex $\Delta(J)$. Finally, we need to verify that the
acuteness condition \eqref{acute} holds for $\mathcal{A}^* =
\{a^*_1,...,a^*_H\}$.

For two $I$-vectors $a$ and $b$, define a weighted taxicab metric 
\begin{equation*}
	\phi_\Omega(a,b) := \sum_{i=1}^I \omega_i |a_i - b_i|,
\end{equation*}
then the standard taxicab metric $\phi$ used above is $\phi_\Omega$ with $%
\Omega = \mathbf{I}_I$. Moreover, letting $a^* = \Omega^{1/2}a$ and $b^* =
\Omega^{1/2}b$, where each of $a$ and $b$ is an $I$-dimensional 0-1 vector,
we have 
\begin{equation*}
	{a^*}^{\prime }b^* = \sum_{i=1}^I \omega_i[1 - |a_i - b_i|] = \bar \omega -
	\phi_\Omega(a,b)
\end{equation*}
with $\bar \omega = \sum_{i=1}^I \omega_i$. Then for every triplet $%
(a_0^*,a_1^*,a_2^*)$ of distinct elements in $\mathcal{A}^*$ 
\begin{eqnarray*}
	(a^*_1 - a^*_0)^{\prime }(a^*_2 - a^*_0) &=& \bar \omega -
	\phi_\Omega(a_1,a_2) - \bar \omega + \phi_\Omega(a_0,a_2) - \bar \omega +
	\phi_\Omega(a_0,a_1) + \bar \omega - \phi_\Omega(a_0,a_0) \\
	&=& \phi_\Omega(a_1,a_2) - \phi_\Omega(a_0,a_2) - \phi_\Omega(a_0,a_1) \\
	&\geq& 0,
\end{eqnarray*}
which is the desired acuteness condition. Since $\mathcal{J}_N$ can be written as the
minimum of the quadratic form with identity-matrix weighting subject to the
cone generated by $a^*$'s, all the previous arguments developed for the case
with $\Omega = \mathbf{I}_I$ remain valid.

Defining $\xi \sim {\mathrm{N}}(0,\hat S)$ and $\zeta = B\xi$, 
\begin{eqnarray*}
	\tilde{\mathcal{J}}_N &\sim& \min_{Bt \leq -\tau_N\phi} N \left[\binom {KB} M
	(\hat \eta_{\tau_N} + N^{-1/2}\xi - t)\right]^{\prime }{P^{-1}}^{\prime
	}\Omega P^{-1}\left[\binom {KB} M (\hat \eta_{\tau_N} + N^{-1/2}\xi - t)%
	\right] 
	\\
	&=& N \min_{\gamma = [{\gamma^\leq}',{\gamma^=}']', \gamma^\leq \leq -\tau_N\bar \phi, \gamma^= = 0, \gamma \in \text{col}(B)} t\left(K\left[%
	\hat \gamma_{\tau_N} + N^{-1/2}\zeta - \gamma\right]\right)
\end{eqnarray*}
conditional on data $\{\{d_{i|j,n}\}_{i=1}^{I_{j}}%
\}_{n=1}^{N_{j}}$, $j=1,2,...,J$.  Moreover, defining $\gamma^\tau = \gamma + \tau_N \phi$ in the above, and
using the definitions of $\beta_{\tau}(\cdot)$ and $s(\cdot)$ 
\begin{eqnarray*}
	\tilde{\mathcal{J}}_N &\sim& N 
	\min_{\gamma^\tau = [{{\gamma^\tau}^\leq}',{{\gamma^\tau}^=}']', {\gamma^\tau}^\leq \leq 0, {\gamma^\tau}^= = 0, {\gamma^\tau} \in \text{col}(B)}
	t\left(K\left[\hat \gamma_{\tau_N} + \tau_N \phi + N^{-1/2}\zeta -
	\gamma^\tau\right]\right) \\
	&=& N 
	\min_{\gamma^\tau = [{{\gamma^\tau}^\leq}',{{\gamma^\tau}^=}']', {\gamma^\tau}^\leq \leq 0, {\gamma^\tau}^= = 0, {\gamma^\tau} \in \text{col}(B)}
	t\left(K\left[%
	\gamma_{\tau_N}(\hat g) + \tau_N \phi + N^{-1/2}\zeta - \gamma^\tau\right]%
	\right) \\
	&=& N 
	\min_{\gamma^\tau = [{{\gamma^\tau}^\leq}',{{\gamma^\tau}^=}']', {\gamma^\tau}^\leq \leq 0, {\gamma^\tau}^= = 0, {\gamma^\tau} \in \text{col}(B)}
	t\left(K\left[%
	\beta_{\tau_N}(\hat g) + N^{-1/2}\zeta - \gamma^\tau\right]\right) \\
	&=& s\left(N^{1/2}\beta_{\tau_N}(\hat g) + \zeta \right)
\end{eqnarray*}
Let $\varphi_N(\xi) := N^{1/2}\beta_{\tau_N}(\tau_N \xi)$ for $\xi =
(\xi_1,...,\xi_m)^{\prime } \in \text{col}(B)$, then from the property of $%
\beta_\tau$ shown above, its $k$-th element $\varphi_N^k$ for $k \leq \bar m$ satisfies 
\begin{equation*}
	\varphi_N^k(\xi) = 0
\end{equation*}
if $|\xi^k| \leq \delta$ and $\xi^j \leq \delta, 1 \leq j \leq m$ for large
enough $N$.  Note $\varphi_N^k(\xi) = N^{1/2}\beta_N^k(\tau_N\xi) =  0$ for $k > \bar m$.   Define $\hat \xi := \hat g / \tau_N$ and using the definition of 
$\varphi_N$, we write 
\begin{equation}  \label{Jins}
	\tilde{\mathcal{J}}_N \sim s\left(\varphi_{N}(\hat \xi) + \zeta \right).
\end{equation}
Now we invoke Theorem 1 of Andrews and Soares (2010, AS henceforth). As noted before, the function $t$ is a
positive definite quadratic form on $\mathbf{R}^\ell$, and so is its
restriction on col$(B)$. Then their Assumptions 1-3 hold for the function $%
s$ defined above if signs are adjusted appropriately as our formulae deal
with negativity constraints, whereas AS work with positivity
constraints. (Note that Assumption 1(b) does not apply here since we use a
fixed weighting matrix.) The function $\varphi_N$ in \eqref{Jins} satisfies
the properties of $\varphi$ in AS used in their proof of Theorem 1. AS
imposes a set of restrictions on the parameter space (see their Equation
(2.2) on page 124). Their condition (2.2) (vii) is a Lyapounov condition for
a triangular array CLT. Following AS, consider a sequence of distributions 
$\pi_N = [\pi_{1 N}^{\prime },...,\pi_{J N}^{\prime }]^{\prime }, N =
1,2,... $ in $\mathcal{P }\cap \mathcal{C}$ such that (1) $\sqrt N B \pi_N \rightarrow
h$ for a non-positive $h$ as $N \rightarrow \infty$ and (2) Cov$%
_{\pi_N}(\sqrt N B\hat \pi) \rightarrow \Sigma$ as $N \rightarrow \infty$
where $\Sigma$ is positive semidefinite.  The Lyapounov condition holds for  $b_k' \hat \pi$  under $\pi_N$ for $k \in \mathcal K^R$ as Condition \ref{condition 1} is imposed for $\pi_N \in \mathcal{P}$.  We do not impose Condition \ref{condition 1} for $k \in \mathcal K^D$.
Note, however, that: (i) The equality $b_k'\hat \pi \leq 0$ holds by construction for every  $k \in \mathcal K^D$ and therefore its behavior does not affect $\mathcal{J}_N$;
(ii) If $\mathrm{var}_{\pi_{N}}(g_k)$  converges to zero for some $k \in \mathcal K^D$, then $\sqrt N b_k'[\tilde \eta_{\tau_N} - \hat \eta_{\tau_N}] = o_p(1)$ and therefore its contribution to $\tilde{\mathcal{J}}_N$ is asymptotically negligible in the size calculation. The other conditions in AS10, namely (2.2)(i)-(vi), hold trivially. Finally, Assumptions GMS 2 and GMS 4 of AS10 are concerned with their thresholding parameter $\kappa_N$ for the $k$-th moment inequality, and by letting $\kappa_N = N^{1/2}\tau_N\phi_k$, the former holds by the condition $\sqrt N \tau_N \uparrow \infty$ and the latter by $\tau_N \downarrow 0$. Therefore we conclude 
\begin{equation*}  \label{eq:size}
	\liminf_{N \rightarrow \infty} \inf_{\pi \in \mathcal{P} \cap \mathcal{C}} \Pr\{\mathcal{J}_N
	\leq \hat c_{1 - \alpha}\} = 1 - \alpha.
\end{equation*}
\end{proof}

\noindent {\bf Further details of the procedure  in Section \ref{sec:extending}} \:  The setting in this section is as follows.  Let $\tilde p_j \in {\bf R}^K_{++}$ denote the unnormalized price vector, fixed for each period $j$.   Let $(S,\mathcal S,P)$ denote the underlying probability space.  Since we have repeated cross-sections over $J$ periods, write $P = \otimes_{j=1}^J P^{(j)}$, a $J$-fold product measure.

We first develop a smoothing procedure based on a series estimator (see, for example, \citeasnoun{newey1997}) for $\pi$ to deal with a situation where total expenditure $W$ is continuously distributed, yet exogenous.  We need some notation and definitions to formally state the asymptotic theory behind our procedure with smoothing.  With exogeneity  we have 
\begin{eqnarray*}
	p_{i|j}(\underline{w}_j)
	&=& \Pr\{D(\tilde p_j/w_{n(j)},u) \in x_{i|j}|w_{n(j)}=\underline{w}_j\}
	\\
	&=& \Pr\{D(\tilde p_j/\underline{w}_j,u) \in x_{i|j}, u \sim P_u\}
\end{eqnarray*}
where the second equality follows from the exogeneity assumption.  Then $\pi_{i|j} =	p_{i|j}(\underline{w}_j)$, which is the estimand in what follows.
Define  $q^K(w) = (q_{1K}(w),...,q_{KK}(w))'$, where $q_{jK}(w), j = 1,...,K$ are basis functions (e.g. power series or splines) of $w$. Instead of sample frequency estimators, for each $j, 1 \leq j \leq J$ we use 
\begin{eqnarray*}
	\hat{\pi}_{i|j} &=& q^{K(j)}(\underline{w}_j)'\widehat Q^{-}(j) \sum_{n(j)=1}^{N_j}q^{K(j)}(w_{n(j)})d_{i|j,n(j)}/N_j,
	\\
	\widehat Q(j) &=& \sum_{n(j)=1}^{N_j} q^{K(j)}(w_{n(j)})q^{K(j)}(w_{n(j)})'/N_j
	\\
	\hat{\pi}_{j} &=&(\hat{\pi}_{1|j},...,\hat{\pi}_{I_{j}|j})^{\prime }, \\
	\hat{\pi} &=&(\hat{\pi}_{1}^{\prime },...,{\hat{\pi}_{J}}^{\prime })^{\prime
	},
\end{eqnarray*}%
to estimate $\pi_{i|j}$, where $A^{-}$ denotes a symmetric generalized inverse of $A$ and $K(j)$ is the number of basis functions applied to Budget $\mathcal{B}_{j}$.   
The estimators $\hat{\pi}_{i|j}$'s may not take their values in $[0,1]$.  This does not seem to cause a problem asymptotically, though as in \citeasnoun{imbens2009}, we may (and do, in the application) instead use 
$$
\hat{\pi}_{i|j}=G\left(q^{K(j)}(\underline{w}_j)'\widehat Q^{-}(j) \sum_{n(j)=1}^{N_j}q^{K(j)}(w_{n(j)})d_{i|j,n(j)}/N_j\right),
$$
where $G$ denotes the CDF of Unif$(0,1)$.    Then an
appropriate choice of $\tau _{N}$ is $\tau_N = \sqrt{\frac{\log \underline n}{\underline n}}$ with 
$$
\underline n = \min_j N_j I_j / {\mathrm{trace}}({v}_N^{(j)})
$$ 
where ${v}_N^{(j)}$ is defined below.  Strictly speaking, asymptotics
with nonparametric smoothing involve bias, and the bootstrap does not solve
the problem. A standard procedure is to claim that one used undersmoothing
and can hence ignore the bias, and we follow this convention.
The bootstrapped test statistic  $%
\tilde{J}_{N}$ is obtained applying the same replacements to the
formula \eqref{J-tilde}, although generating $\tilde{\eta}_{\tau _{N}}$
requires a slight modification. Let $\hat{\eta}_{\tau _{N}}(j)$ be the $j$%
-th block of the vector $\hat{\eta}_{\tau _{N}}$, and ${\hat v}_N^{(j)}$ satisfy 
${\hat v}_N^{(j)}{v_{N}^{(j)}}^{-1} \rightarrow_p {\bf I}_{I_j}$, where
$$
v_{N}^{(j)} = [{\bf I}_{I_j} \otimes q^{K(j)}(\underline{w}_j)'Q_{N}(j)^{-1}]\Lambda_{N}^{(j)}[{\bf I}_{I_j} \otimes Q_{N}^{-1}(j)q^{K(j)}(\underline{w}_j)]
$$
with $Q_{N}(j) := \mathrm E[q^{K(j)}(w_{n(j)}) {q^{K(j)}(w_{n(j)})}']$, $\Lambda_{N}^{(j)} := \mathrm E[\Sigma^{(j)}(w_{n(j)}) \otimes q^{K(j)}(w_{n(j)}) {q^{K(j)}(w_{n(j)})}']$, and  $\Sigma^{(j)}(w) := \mathrm{Cov}[d_{j,n(j)}|w_{n(j)} = w] $.  Note that $\Sigma^{(j)}(w) = \mathrm{diag}\left(p^{(j)}(w)\right) -
p^{(j)}(w) p^{(j)}(w)^{\prime }$ where  $p^{(j)}(w) = [p_{1|j}(w),...,p_{I_j|j}(w)]'$.  For example, one may use 
$$
\hat v_{N}^{(j)} = [{\bf I}_{I_j} \otimes q^{K(j)}(\underline{w}_j)'\hat Q^{-}(j)]\widehat \Lambda(j)[{\bf I}_{I_j} \otimes \hat Q^{-}(j)q^{K(j)}(\underline{w}_j)]
$$
with $\widehat \Lambda(j) = \frac 1 {N_j}\sum_{n(j) = 1}^{N_j} \left[\widehat \Sigma^{(j)}(w_{n(j)}) \otimes q^{K(j)}(w_{n(j)}) {q^{K(j)}(w_{n(j)})}'\right]$, $\widehat \Sigma^{(j)}(w) = \mathrm{diag}\left(\widehat p^{(j)}(w)\right) -
\widehat p^{(j)}(w) \widehat p^{(j)}(w)^{\prime }$, $\widehat p^{(j)}(w) = [\widehat p_{1|j}(w),...,\widehat p_{I_j|j}(w)]'$ and 
$\widehat p_{i|j}(w) = q^{K(j)}(w)'\widehat Q^{-}(j) \sum_{n(j)=1}^{N_j}q^{K(j)}(w_{n(j)})d_{i|j,n(j)}/N_j$.  We use $\tilde{\eta}_{\tau _{N}}=(\tilde{\eta}%
_{\tau _{N}}(1)^{\prime },...,\tilde{\eta}_{\tau _{N}}(J)^{\prime })'$ for
the smoothed version of $\tilde{J}_{N}$, where $\tilde{\eta}%
_{\tau _{N}}(j):=\hat{\eta}_{\tau _{N}}(j)+\frac{1}{\sqrt{N_{j}}}N(0,{\hat v}_N^{(j)}),j=1,...,J$.

Noting  $\{\{d_{i|j,n(j)}\}_{i=1}^{I_j},w_{n(j)}\}_{n(j)=1}^{N_j}$   are IID-distributed within each time period $j, 1 \leq j \leq J$, let $(d_j,w_j)$ denote the choice-log-expenditure pair of a consumer facing budget $j$.     Let $d = [d_1',...,d_J']'$ and  ${\bf w} = [w_1,...,w_J]'$, and define $g = Bd = [g_1,...,g_m]'$ as before.
Let $\mathcal W_j$ denote the support of $w_{n(j)}$.  For a symmetric matrix $A$, $\lambda_{\mathrm{min}}$ signifies its smallest eigenvalue.         
\begin{condition}
	\label{condition:kernel}  There exist positive constants $c_1$,  $c_2$,  $\delta$, and $\zeta(K)$,  $K \in \mathbf N$ such that 
	the following holds:
	
	\begin{itemize}
		
		\item[{(i)}] $\pi \in \mathcal C$;
		
		\item[{(ii)}] For each  $k \in {\mathcal K}^R$,   {\rm{var}}$(g_k|{\bf w} = (\underline{w_1},...,\underline{w_J})')  \geq s^2(F_1,...,F_J)$ and $\mathrm{E}[(g_k/s(F_1,...,F_J))^4|{\bf w} = (\underline{w_1},...,\underline{w_J})'] < c_1$ hold for every $(\underline{w_1},...,\underline{w_J}) \in {\mathcal W}_1 \times \cdots  {\mathcal W}_J $; 
		
		\item[{(iii)}] $
		\sup_{w \in \mathcal W_j}|p_{i|j}(w) - q^K(w)'\beta_{K}^{(j)}|   \leq c_1 K^{-\delta}  
		$
		holds with some $K$-vector $\beta_{K}^{(j)}$ for every $K \in \bf N$, $ 1 \leq i \leq I_j, 1 \leq j \leq J$;
		
		\item[{(iv)}] Letting $\widetilde q^K := C_{K,j} q^K$, $\lambda_{\mathrm{min}}\mathrm E[\widetilde q^{K}(w_{n(j)}) \widetilde q^{K}(w_{n(j)})'] \geq c_2$ holds for every $K$ and $j$, where $C_{K,j}, K \in {\bf N}, 1 \leq j \leq J$ are constant nonsingular matrices;
		
		\item[{(v)}] $\max_j \sup_{w \in \mathcal W_j}\|\widetilde q^K(w)\| \leq c_2 \zeta(K)$ for every $K \in \bf N$.
		
	\end{itemize}
	
\end{condition}

\qed

\noindent Condition \ref{condition:kernel}(ii) is a version of Condition \ref{condition 1} that accommodates the conditioning by $w$ and series estimation.  Conditions \ref{condition:kernel}(iii)-(v) are standard regularity commonly used in the series regression literature: (iii) imposes a uniform approximation error bound, (iv) avoids singular design (note the existence of the matrices $C_{K,j}$ suffices) and (v) controls the lengths of the series terms used.

The next condition imposes restrictions on tuning parameters.

\begin{condition}
\label{cond:tuning1}	
$\tau_N$ and $K(j), j = 1,...,J$ satisfy $\sqrt{N_j} K^{-\delta}(j) \downarrow 0$, $\zeta(K(j))^2K(j)/N_j \downarrow 0$, $j = 1,...,J$, $\tau_N
\downarrow 0$, and $\sqrt{\underline n}
\tau_N \uparrow \infty$.

\end{condition}

\begin{proof}[\textup{\textbf{Proof of Theorem~\protect\ref{thmSM}}}]
We begin by introducing some notation.
\begin{notation}  Let $b_{k,i}$, $k = 1,...,m$, $i = 1,...,I$ denote the $(k,i)$ element of $B$%
	, then define 
	\begin{equation*}
	b_{k}(j) = [b_{k, {N_1 + \cdots N_{j-1} + 1}}, b_{k, {N_1 + \cdots N_{j-1} +
			2}}, ... , b_{k, {N_1 + \cdots N_{j}}}]^{\prime }
	\end{equation*}
	for $1 \leq j \leq J$ and $1 \leq k \leq m$.
	Let $B^{(j)} : = [b_1(j),...,b_m(j)]^{\prime }\in \mathbf{R}^{m \times {I_j}%
	} $. For $F \in \mathcal{F}$ and $1 \leq j \leq J$, define 
	\begin{equation*}
		p_{F}^{(j)}(w) := E_{F}[d_{j,n(j)}|w_{n(j)} = w], \quad \pi_{F}^{(j)} =
		p_{F}^{(j)}(\underline{w}_j), \quad \quad \pi_F = [{\pi_{F}^{(1)}}^{\prime
		},...,{\pi_{F}^{(J)}}^{\prime }]^{\prime }
	\end{equation*}
	and 
	\begin{equation*}
		\Sigma_{F}^{(j)}(w) := \mathrm{Cov}_{F}[d_{j,n(j)}|w_{n(j)} = w].
	\end{equation*}
\end{notation}
Note that $\Sigma_{F}^{(j)}(w) = \mathrm{diag}\left(p_F^{(j)}(w)\right) -
p_F^{(j)}(w) p_F^{(j)}(w)^{\prime }. $ 

The proof mimics the proof of Theorem \ref{thm1}, except for the treatment of $\hat \pi$. Instead of the sequence 
	$\pi_N, N = 1,2,...$ in $\mathcal{P }\cap \mathcal C$, consider a sequence of
	distributions $F_N = [F_{1 N},...,F_{J N}], N = 1,2,...$ in $\mathcal{F}$
	such that  $\sqrt{N_j/K(j)}B^{(j)}\pi_{F_N}^{(j)} \rightarrow h_j, h_j \leq 0, 1 \leq j \leq J$ as $N \rightarrow \infty$.  Define $Q_{F_N}^{(j)} = \mathrm E_{F_N}[q^{K(j)}(w_{n(j)}) {q^{K(j)}}(w_{n(j)})']$ and $\Xi_{F_N}^{(j)} = \mathrm E_{F_N}[B^{(j)}\Sigma_{F_N}^{(j)}(w_{n(j)}){B^{(j)}}' \otimes q^{K(j)}(w_{n(j)}) {q^{K(j)}}(w_{n(j)})']$, and let 
	$$
	V_{F_N}^{(j)} := [{\bf I}_{m} \otimes q^{K(j)}(\underline{w}_j)'{Q_{F_N}^{(j)}}^{-1} ]\Xi_{F_N}^{(j)}[{\bf I}_{m} \otimes {Q_{F_N}^{(j)}}^{-1}q^{K(j)}(\underline{w}_j)]
	$$
	and 
	$$
	V_{F_N} := \sum_{j=1}^J V_{F_N}^{(j)}. 
	$$
	Then by adapting the proof of Theorem 2 in \citeasnoun{newey1997} to the triangle array for the repeated crosssection setting, we obtain 
	$$
	\sqrt N {V_{F_N}}^{- \frac 1 2}B[\hat \pi - \pi_{F_N}] \overset{F_N}{\leadsto}
	N(0, {\bf I}_{m}).
	$$
	The rest is the same as the proof of Theorem \ref{thm1}.
\end{proof}

\
\

Next we turn to the definition of our endogeneity corrected estimator $\widehat{\pi_{\rm{EC}}}$ propose a bootstrap algorithm for it, and show its validity.  Exogeneity of budget sets is a standard assumption in classical demand analysis based on random utility models; for example, it is assumed, at least implicitly, in \citeasnoun{mcfadden-richter}. Nonetheless, the assumption can be a concern in applying our testing procedure to a data set such as ours. Recall that the budget sets $\{\mathcal{B}_j\}_{j=1}^J$ are based on prices and total expenditure. The latter is likely to be endogenous, which should be a concern to the econometrician.         

As independence between utility and budgets is fundamental to McFadden-Richter theory, addressing it in our testing procedure might seem difficult. Fortunately, recent advances in nonparametric identification and estimation of models with endogeneity inform a solution.  To see this, it is useful to rewrite the model so that we can cast it into a framework of nonseparable models with endogenous covariates. Writing $p_j = \tilde{p}_j/W$, where $\tilde{p}_j$ is the unnormalized price vector, the essence of the problem is as follows: Stochastic rationalizability imposes restrictions on the distributions of $y=D(p,u)$ for different $p$ when $u$ is distributed according to its population marginal distribution $P_u$, but the observed conditional distribution of $y$ given $p$ does not estimate this when $w$ and $u$ are interrelated.   
In particular, if we define $\mathcal J_{\rm{EC}}  = \min_{\nu \in {\bf R}^h_+}[\pi_{\rm{EC}}  - A\nu]'\Omega[\pi_{\rm{EC}}  - A\nu]$,  with the definition of $\pi_{\rm{EC}}$ provided in Section \ref{sec:extending}, then $\mathcal J_{\rm{EC}}=0$ iff stochastic rationalizability holds. Note that the new definition $\pi_{\rm{EC}}$ recovers the previous definition of $\pi$  when $w$ is exogenous.

Our estimator uses the control function approach.  For example, given a reduced form $w = h_j(z,\varepsilon)$ with $h_j$ monotone in $\varepsilon$ and $z$ is an instrument, one may use 
\begin{equation}\label{eq:cfexample}
\varepsilon = F_{w|z}^{(j)}(w|z)
\end{equation}
where $F_{w|z}^{(j)}$ denotes the conditional CDF of $w$ given $z$ under $P^{(j)}$ when the random vector $(w,z)$ obeys the probability law $P^{(j)}$; see \citeasnoun{imbens2009} for this type of control variable in the context of cross-sectional data.   
Note that $\varepsilon \sim \mathrm{Uni}(0,1)$ under every $P^{(j)}, 1 \leq j \leq J$ by construction.    Let $P^{(j)}_{y|w,\varepsilon}$ denote the conditional probability measure for $y$ given $(w,\varepsilon)$ corresponding to $P^{(j)}$.  Adapting the argument in \citeasnoun{imbens2009} and \citeasnoun{blundell2003}, under the assumption that supp$(w) = $ supp$(w| \varepsilon)$ under $P^{(j)}, 1 \leq j \leq J$ we have        
\begin{eqnarray*}
	\pi(p_j,x_{i|j}) 
	&=& \int_0^1 \int_u {\bf 1}\{D_j(\underline{w}_j,u) \in x_{i|j} \} dP_{u|\varepsilon}^{(j)}d\varepsilon 
	\\
	&=& \int_0^1 {P^{(j)}_{y|w,\varepsilon}}\left\{y \in x_{i|j}|w = \underline{w}_j,\varepsilon \right\}d\varepsilon, \quad 1 \leq j \leq J.
\end{eqnarray*}
This means that $\pi_{\rm{EC}}$ can be estimated nonparametrically.  

 To estimate $\widehat{\pi_{\mathrm{EC}}}$, we can proceed in two steps as follows.  
The first step is to obtain control variable estimates $\widehat \epsilon_{n(j)}, n(j) = 1,...,N_j$ for each  $j$.  For example, let 
$
\widehat F_{w|z}^{(j)}
$
be a nonparametric estimator for $F_{w|z}$ for a given instrumental variable $z$ in period $j$. For concreteness, we consider a series estimator as in \citeasnoun{imbens-newey}.    Let $r^L(z) = (r_{1L}(z),...,r_{LL}(z))$, where $r_{\ell L}(z), \ell = 1,...,L$ are basis functions, then define 
$$
\widehat F_{w|z}^{(j)}(w|z) = r^{L}(z)'\widehat R^-(j)\sum_{n(j)=1}^{N_j}r^{L(j)}(z_{n(j)}){\bf 1}\{w_{n(j)} \leq w\}/N_j
$$
where
$$
\widehat R(j) = \sum_{n(j)=1}^{N_j} r^{L(j)}(z_{n(j)})r^{L(j)}(z_{n(j)})'/N_j.
$$
Let
$$
\widetilde \epsilon_{n(j)} = \widehat F_{w|z}^{(j)}(w_{n(j)}|z_{n(j)}), n(j) = 1,...,N_j.
$$
Choose a sequence $\upsilon_N \rightarrow 0, \upsilon_N > 0$ and define $\iota_N(\varepsilon) = (\varepsilon + \upsilon_N)^2/4\upsilon_N$, then let   
$$
\gamma_N(\varepsilon) = 
\begin{cases}
1 & {\text{if } } \varepsilon > 1 + \upsilon_N 
\\
1 - \iota_N(1 - \varepsilon) & {\text{if } } 1 - \upsilon_N < \varepsilon \leq 1 + \upsilon_N
\\
\varepsilon  & {\text{if } } \upsilon_N \leq \varepsilon \leq 1 - \upsilon_N
\\
\iota_N(\varepsilon) & {\text{if } }   - \upsilon_N \leq \varepsilon \leq \upsilon_N
\\
0 & {\text{if } } \varepsilon < -\upsilon_N
\end{cases}
$$
then our control variable is 
$
\widehat \varepsilon_{n(j)} = \gamma_N(\widetilde{\varepsilon}_{n(j)}), n(j) = 1,...,N_j.
$

The second step is nonparametric estimation of ${P^{(j)}_{y|w,\varepsilon}}\left\{y \in x_{i|j}|w = \underline{w}_j,\varepsilon \right\}$.   Let $\widehat \chi_{n(j)} = (w_{n(j)},\widehat \varepsilon_{n(j)})'$, $n(j) = 1,...,N_j$ for each $j$.  Write $s^{M(j)}(\chi) = (s_{1M(j)}(\chi),...,s_{M(j)M(j)}(\chi))'$, where $s_{m M(j)}(\chi), \chi \in {\bf R}^{K+1}, m = 1,...,M(j)$ are basis functions, then our estimator for ${P^{(j)}_{y|w,\varepsilon}}\left\{y \in x_{i|j}|w = \cdot,\varepsilon = \cdot \right\}$ evaluated at $\chi = (w,\varepsilon)$ is
\begin{eqnarray*}
	\widehat  {P^{(j)}_{y|w,\varepsilon}}\left\{y \in x_{i|j}|w,\varepsilon \right\}         &=& s^{M(j)}(\chi)'\widehat S^-(j)\sum_{n(j)=1}^{N_j}s^{M(j)}(\widehat \chi_{n(j)})d_{i|j,n(j)}/N_j     
	\\
	&=& s^{M(j)}(\chi)'\widehat{\alpha}_i^{M(j)} 
\end{eqnarray*}
where
$$
\widehat S(j) = \sum_{n(j)=1}^{N_j} s^{M(j)}(\widehat \chi_{n(j)})s^{M(j)}(\widehat \chi_{n(j)})'/N_j, \quad  \widehat{\alpha}_i^{M(j)} := \widehat S^-(j)\sum_{n(j)=1}^{N_j}s^{M(j)}(\widehat \chi_{n(j)})d_{i|j,n(j)}/N_j.
$$
Our endogeneity corrected conditional probability $\pi(p_j,x_{i|j})$ is a linear functional of ${P^{(j)}_{y|w,\varepsilon}}\left\{y \in x_{i|j}|w = \underline{w}_j,\varepsilon \right\}$, thus plugging in $\widehat {P^{(j)}_{y|w,\varepsilon}}\left\{y \in x_{i|j}|w = \underline{w}_j,\varepsilon \right\}$ into the functional, we define
\begin{eqnarray*}
	\widehat{\pi(p_j,x_{i|j})} &:=&
	\int_0^1  \widehat {P^{(j)}_{y|w,\varepsilon}}\left\{y \in x_{i|j}|w = 
	\underline{w}_j,\varepsilon \right\} d\varepsilon 
	\\
	&=& D(j)'\widehat{\alpha}_i^{M(j)},  
	\\
	\text{where } \quad D(j) &:=& \int_0^1 s^{M(j)}\left(
	\begin{bmatrix}
		\underline{w}_j \\
		\varepsilon
	\end{bmatrix}
	\right)d\varepsilon \qquad \qquad \qquad i = 1,...,I_j, j = 1,...,J
\end{eqnarray*}
and 
$$
\widehat{\pi_{\rm{EC}}} = [\widehat{\pi(p_1,x_{1|1})},...,\widehat{\pi(p_{1},x_{I_1|1})},\widehat{\pi(p_2,x_{1|2})},...,\widehat{\pi(p_{2},x_{I_2|2})},...,\widehat{\pi(p_J,x_{1|J})},...,\widehat{\pi(p_{J},x_{I_J|J})}]'.
$$ 
The final form of the test statistic  is
$$
\mathcal J_{{\rm EC}_N} = N \min_{\nu \in {\bf R}^h_+}[\widehat{\pi_{\rm{EC}}} - A\nu]'\Lambda[\widehat{\pi_{\rm{EC}}} - A\nu].
$$
The calculation of critical values can be carried out in the same way as the testing procedure with the series estimator $\hat \pi$ for the exogenous case, though the covariance matrix $v_N^{(j)}$ needs modification.  With the nonparametric endogeneity correction, the modified version of $v_N^{(j)}$ is
$$
\bar v_N^{(j)} = [{\bf I}_{I_j} \otimes D(j)'S_N(j)^{-1}] \bar \Lambda_N^{(j)} [{\bf I}_{I_j} \otimes S_N(j)^{-1}D(j)]
$$
where 
$$
S_{N}(j) = \mathrm E[s^{M(j)}(\chi_{n(j)}) {s^{M(j)}(\chi_{n(j)})}'],
\quad
\bar \Lambda_N^{(j)} = \bar \Lambda_{1_N}^{(j)} + \bar \Lambda_{2_N}^{(j)},
$$ 
$$
\bar \Lambda_{1_N}^{(j)} = {\mathrm{E}}[\bar \Sigma^{(j)}(\chi_{n(j)}) \otimes s^{M(j)}(\chi_{n(j)})s^{M(j)}(\chi_{n(j)})'], \quad \bar \Lambda_{2_N}^{(j)} =  {\mathrm{E}}[m_{n(j)}m_{n(j)}']
$$ 
with 
$$
\bar \Sigma^{(j)}(\chi) := {\mathrm{Cov}}[d_{j,n(j)}|\chi_{n(j)} = \chi],
$$
$$
m_{n(j)} := [m_{1,n(j)}', m_{2,n(j)}',\cdots,m_{I_j,n(j)}']',
$$
\begin{eqnarray*}
	&&m_{i,n(j)} :=
	\\
	&&{\mathrm E}\left[\dot \gamma_N(\varepsilon_{m(j)}) \frac{\partial}{\partial \varepsilon}{P^{(j)}_{y|w,\varepsilon}}\left\{y \in x_{i|j}|w_{m(j)},\varepsilon_{m(j)} \right\}s^{M(j)}(\chi_{m(j)}) r^{L(j)}(z_{m(j)})'R_N(j)^{-1}r^{L(j)}(z_{n(j)})u_{mn(j)}
	\phantom{|d_{i|j,n(j)}, w_{n(j)}, z_{n(j)}  } 
	\right] 
	\\
	&& \qquad
	\left|
	d_{i|j,n(j)}, w_{n(j)}, z_{n(j)}  
	\phantom{\frac{\partial}{\partial \varepsilon}{P^{(j)}_{y|w,\varepsilon}}}
	\!\!\!\!\!\!\!\!\!\!\!\!\!\!\!\!\!\!\!\!   \right], 
\end{eqnarray*}
$$
R_N(j) := {\mathrm E}[r^{L(j)}(z_{n(j)})r^{L(j)}(z_{n(j)})'], \quad u_{mn(j)} := {\bf 1}\{w_{n(j)} \leq w_{m(j)}\} - F_{w|z}^{(j)}(w_{m(j)}|z_{n(j)}). 
$$
Define
$$
\underline n_{\mathrm{EC}} = \min_j N_j I_j / {\mathrm{trace}}(\bar{v}_N^{(j)}),
$$ 
then a possible choice for $\tau_N$ is $\tau_N = \sqrt{\frac{\log \underline n_{\mathrm{EC}}}{\underline n_{\mathrm{EC}}}}$.  Proceed as for $\tilde {\mathcal J}_N$ earlier in this section, replacing $\hat v_N^{(j)}$ with a consistent estimator for $\bar v_N^{(j)}$ for $j = 1,...,J$, to define the bootstrap version $\tilde {\mathcal J}_{\rm EC}$.

We  impose some conditions to show the validity of the endogeneity-corrected test.  Let $\varepsilon_{n(j)}$ be the value of the control variable $\varepsilon$ for the $n(j)$-th consumer facing budget $j$.  Noting  $\{\{d_{i|j,n(j)}\}_{i=1}^{I_j},w_{n(j)},\varepsilon_{n(j)}  \}_{n(j)=1}^{N_j}$   are IID-distributed within each time period $j, 1 \leq j \leq J$, let $(d_j,w_j,\varepsilon_j)$ denote the choice-log expenditure-control variable triplet of a consumer facing budget $j$.     Let $d = [d_1',...,d_J']'$, ${\bf w} = [w_1,...,w_J]'$, ${\bf e} = [\varepsilon_1,...,\varepsilon_J]'$ and define $g = Bd = [g_1,...,g_m]'$ as before.  Let $\mathcal X_j = \mathrm{supp}(\chi_{n(j)})$, $\mathcal Z_j = \mathrm{supp}(z_{n(j)})$, and $\mathcal E_j = \mathrm{supp}(\varepsilon_{n(j)})$, $1 \leq j \leq J$.  Following the above discussion, define an ${\bf R}^{I}$-valued functional 
$$
\pi(P_{y|w,\varepsilon}^{(1)},...,P_{y|w,\varepsilon}^{(J)}) = [\pi_{1|1}(P_{y|w,\varepsilon}^{(1)}),...,\pi_{I_1|1}(P_{y|w,\varepsilon}^{(1)}),  \pi_{1|2}(P_{y|w,\varepsilon}^{(2)}),...,\pi_{I_2|2}(P_{y|w,\varepsilon}^{(2)}),
...,
\pi_{1|J}(P_{y|w,\varepsilon}^{(J)}),...,\pi_{I_J|J}(P_{y|w,\varepsilon}^{(J)})]'
$$
where 
$$
\pi_{i|j}(P_{y|w,\varepsilon}^{(j)}) := \int_0^1 {P^{(j)}_{y|w,\varepsilon}}\left\{y \in x_{i|j}|w = \underline{w}_j,\varepsilon \right\}d\varepsilon 
$$
and $\varepsilon_{n(j)} := F_{w|z}^{(j)}(w_{n(j)}|z_{n(j)})$ for every $j$.  

\begin{condition}
	\label{condition:EC}  There exist positive constants $c_1$, $c_2$, $\delta_1$, $\delta$, $\zeta_r(L)$, $\zeta_s(M)$, and $\zeta_1(M)$,  $L \in \mathbf N$, $M \in \mathbf N$ such that 
	the following holds:
	
	\begin{itemize}
		
		\item[{(i)}] The distribution of $w_{n(j)}$ conditional on $z_{n(j)} = z$ is continuous for every $z \in \mathcal Z_j$, $1 \leq j \leq J$;

		\item [(ii)] $\mathrm{supp}(w_{n(j)}|\varepsilon_{n(j)} = \varepsilon) = \mathrm{supp}(w_{n(j)})$ for every $\varepsilon \in [0,1]$, $1 \leq j \leq J$;
		
		\item[{(iii)}] $\pi(P_{y|w,\varepsilon}^{(1)},...,P_{y|w,\varepsilon}^{(J)}) \in \mathcal C$;

		\item[{(iv)}]  For each  $k \in {\mathcal K}^R$, {\rm{var}}$(g_k|{\bf w} = (\underline{w_1},...,\underline{w_J})', {\bf e} = (\underline{\varepsilon_1},...,\underline{\varepsilon_J}))'  \geq s^2(F_1,...,F_J)$ and $\mathrm{E}[(g_k/s(F_1,...,F_J))^4$ $|{\bf w} = (\underline{w_1},...,\underline{w_J})', {\bf e} = (\underline{\varepsilon_1},...,\underline{\varepsilon_J})'] < c_1$ hold for every $(\underline{w_1},...,\underline{w_J},\underline{\varepsilon_1},...,\underline{\varepsilon_J}) \in {\mathcal W}_1 \times \cdots  {\mathcal W}_J \times {\mathcal E}_1 \times \cdots  {\mathcal E}_J$;

		\item[{(v)}] Letting $\widetilde{r}^L := C_{L,j} r^L$, $\lambda_{\mathrm{min}}\mathrm E[\widetilde r^{L}(z_{n(j)})\widetilde r^{L}(z_{n(j)})] \geq c_2$ holds for every $L$ and $j$, where $C_{L,j}, L \in \bf N$, $1 \leq j \leq J$, are constant nonsingular matrices;
		
		\item[{(vi)}] $\max_j\sup_{z \in \mathcal Z_j}\|\widetilde{r}^L(z)\| \leq c_1 \zeta_r(L)$ for every $L \in \bf N$.	
		
		\item[{(vii)}]  $\sup_{w \in \mathcal W_j, z \in \mathcal Z_j}|F_{w|z}^{(j)}(w,z)  - r^L(z)'\alpha_{L}^{(j)}(w)|   \leq c_1 L^{-\delta_1},  1 \leq j \leq J
		$
		holds with some $L$-vector $\alpha_{L}^{(j)}(\cdot)$ for every $L \in \bf N$, $1 \leq j \leq J$;	
		
		\item[{(viii)}] Letting $\widetilde{s}^M := \bar C_{M,j} s^M$, $\lambda_{\mathrm{min}}\mathrm E[\widetilde s^{M}(\chi_{n(j)})\widetilde s^{M}(\chi_{n(j)})] \geq c_2$ holds for every $M$ and $j$, where $\bar C_{M,j}, M \in \bf N$, $1 \leq j \leq J$, are constant nonsingular matrices;
		
		\item[{(ix)}]  $\max_j\sup_{\chi \in \mathcal  X_j}\|\widetilde{s}^M(\chi)\| \leq C \zeta_s(M)$ and $\max_j\sup_{\chi \in \mathcal  X_j}\|\partial\widetilde{s}^M(\chi)/\partial\varepsilon \| \leq c_1 \zeta_1(M)$ and  $\zeta_s(M) \leq C \zeta_1(M)$ for every $M \in \bf N$;
		
		\item[{(x)}]  $\sup_{\chi \in \mathcal X_j} \left|P_{y|w,\varepsilon}^{(j)}\{y \in x_{i|j}|w,\varepsilon \}  - s^M(\chi)'g_{M}^{(i,j)}\right|   \leq c_1 M^{-\delta},  $
		holds with some $M$-vector $g_{M}^{(i,j)}$ for every $M \in \bf N$, $1 \leq i \leq I_j,  1 \leq j \leq J$;

		\item[{(xi)}]  $P_{y|w,\varepsilon}^{(j)}\{y \in x_{i|j}|w,\varepsilon \}$, $1 \leq i \leq I_j,  1 \leq j \leq J$, are twice continuously differentiable in $\chi = (w,\varepsilon)$.  Moreover, 
		$\max_{1 \leq j \leq J}\max_{1 \leq i \leq I_j}\sup_{\chi \in \mathcal \chi_j}\left \|\frac{\partial}{\partial \chi}{P^{(j)}_{y|w,\varepsilon}}\left\{y \in x_{i|j}|w,\varepsilon\right\}\right \| \leq c_1$ and \\
		$\max_{1 \leq j \leq J}\max_{1 \leq i \leq I_j}\sup_{\chi \in \mathcal \chi_j}\left \|\frac{\partial^2}{\partial \chi \partial \chi'}{P^{(j)}_{y|w,\varepsilon}}\left\{y \in x_{i|j}|w,\varepsilon\right\}\right \| \leq c_1$.

	\end{itemize}
	
\end{condition}

\qed

\noindent Since we use the control function approach to deal with potential endogeneity in $w$ (income), Conditions \ref{condition:EC}(i)-(ii) are essential.  See \citeasnoun{blundell2003} and \citeasnoun{imbens2009} for further discussion on these types of restrictions. Just like Condition \ref{condition:kernel}(ii), Condition \ref{condition:EC}(iv) is a version of Condition \ref{condition 1} that accommodates the two-step series estimation.   Conditions \ref{condition:EC}(iv)-(xi) corresponds to standard regularity conditions stated in the context of the two-step approach adopted in this section:    (iv) and (x) imposes uniform approximation error bounds, (v) and (viii) avoid singular designs (note the existence of the matrices $C_{L,j}$ and $\bar C_{M,j}$ suffices), and (vi) and (ix) control the lengths of (the derivatives of) the series terms used.    Condition \ref{condition:EC}(xi)  imposes reasonable smoothness restrictions on the (observable) conditional probabilities  $P_{y|w,\varepsilon}^{(j)}\{y \in x_{i|j}|w,\varepsilon \}$, $1 \leq i \leq I_j,  1 \leq j \leq J$.  

The next condition impose restrictions on tuning parameters.    
\begin{condition}
\label{cond:tuning2}
	Let $\tau_N$, $M(j)$ and $L(j)$, $j = 1,...,J$ satisfy  $\tau_N
	\downarrow 0$, $\sqrt{\underline n_{\mathrm{EC}}}
	\tau_N \uparrow \infty$, $N_j L(j)^{1 - 2\delta_1} \downarrow 0$,  $N_j M(j)^{- 2\delta} \downarrow 0$, $M(j) \zeta_1(M(j))^2 L^2(j)/N_j \downarrow 0$, $\zeta_s(M(j))^6 L^4(j)/N_j \downarrow 0$, and $\zeta_1(M(j))^4 \zeta_r(L(j))^4/N_j \downarrow 0$ and also $\underline{C}(L(j) /N_j + L(j)^{1 - 2 \delta_1}) \leq \upsilon_N^3 \leq \overline{C}(L(j) /N_j + L(j)^{1 - 2 \delta_1})$, for some $0 < \underline{C} < \overline{C}$.  
\end{condition}

\
\

	\begin{proof}[\textup{\textbf{Proof of Theorem~\protect\ref{thmEC}}}] 
		The proof follows the same steps as those in the proof of Theorem \ref{thm1}%
		, except for the treatment of the estimator for $\pi$.   Therefore, instead of the sequence 
		$\pi_N, N = 1,2,...$ in $\mathcal{P }\cap \mathcal C$, consider a sequence of
		distributions $F_N = [F_{1 N},...,F_{J N}], N = 1,2,...$ in $\mathcal{F}_{\rm EC}$ and the corresponding conditional distributions $P_{{y|w,\varepsilon};F_N}^{(j)}\{y \in x_{i|j}|w,\varepsilon \}$ and $F_{{w|z}_N}^{(j)}$, $1 \leq i \leq I_j, 1 \leq j \leq J$, $N = 1,2,...$
		such that  $\sqrt{N_j/(M(j) \vee L(j)))} B^{(j)}\pi_{F_N}^{(j)} \rightarrow h_j, h_j \leq 0, 1 \leq j \leq J$ as $N \rightarrow \infty$, where
		$\pi_{F_N} = \pi(P_{{y|w,\varepsilon}_N}^{(1)},...,P_{{y|w,\varepsilon}_N}^{(J)})$ whereas the definitions of
		$\overline{V}_{F_N}^{(j)}, 1 \leq j \leq J$ are given shortly.  Define $S_{F_N}^{(j)} = \mathrm E_{F_N}[s^{M(j)}(\chi_{n(j)}) {s^{M(j)}}(\chi_{n(j)})']$ as well as  
		$$\bar \Xi_{1_{F_N}}^{(j)} = \mathrm E_{F_N}[B^{(j)}\bar \Sigma_{F_N}^{(j)}(\chi_{n(j)}){B^{(j)}}' \otimes s^{M(j)}(\chi_{n(j)}) {s^{M(j)}}(\chi_{n(j)})']$$
		and 
		$$
		\bar \Xi_{2_{F_N}}^{(j)} =  [B^{(j)} \otimes {\bf I}_{M(j)}]\mathrm{E}_{F_N}[m_{n(j);F_N} m_{n(j);F_N}'][{B^{(j)}}' \otimes {\bf I}_{M(j)}]
		$$ 
		where
		$$
		\Sigma^{(j)}_{F_N}(\chi) := {\mathrm{Cov}_{F_N}}[d_{j,n(j)}|\chi_{n(j)} = \chi],
		$$
		$$
		m_{n(j);F_N} := [m_{1,n(j);F_N}', m_{2,n(j);F_N}',\cdots,m_{I_j,n(j);F_N}']',
		$$
		\begin{eqnarray*}
			&&m_{i,n(j);F_N} :=
			\\
			&&{\mathrm E}_{F_N}\left[\dot \gamma_N(\varepsilon_{m(j)}) \frac{\partial}{\partial \varepsilon}{P^{(j)}_{y|w,\varepsilon;F_N}}\left\{y \in x_{i|j}|w_{m(j)},\varepsilon_{m(j)} \right\}s^{M(j)}(\chi_{m(j)}) r^{L(j)}(z_{m(j)})'R_{F_N}(j)^{-1}r^{L(j)}(z_{n(j)})u_{mn(j);F_N}
			\phantom{|d_{i|j,n(j)}, w_{n(j)}, z_{n(j)}  } 
			\right] 
			\\
			&& \qquad
			\left|
			d_{i|j,n(j)}, w_{n(j)}, z_{n(j)}  
			\phantom{\frac{\partial}{\partial \varepsilon}{P^{(j)}_{y|w,\varepsilon}}}
			\!\!\!\!\!\!\!\!\!\!\!\!\!\!\!\!\!\!\!\!   \right], 
		\end{eqnarray*}
		$$
		R_{F_N}(j) := {\mathrm E_{F_N}}[r^{L(j)}(z_{n(j)})r^{L(j)}(z_{n(j)})'], \quad u_{mn(j);F_N} := {\bf 1}\{w_{n(j)} \leq w_{m(j)}\} - F_{{w|z}_N}^{(j)}(w_{m(j)}|z_{n(j)}). 
		$$
		With these definitions, let 
		$$
		\overline{V}_{F_N}^{(j)} := \left[{\bf I}_{m} \otimes D(j)'{S_{F_N}^{(j)}}^{-1} \right]\bar \Xi_{{F_N}}^{(j)}\left[{\bf I}_{m} \otimes {S_{F_N}^{(j)}}^{-1}D(j)\right]
		$$
		with $\bar \Xi_{{F_N}}^{(j)} = \bar \Xi_{1_{F_N}}^{(j)} + \bar \Xi_{2_{F_N}}^{(j)}$.  Define 
			$$
			\overline{V}_{F_N} := \sum_{j=1}^J \overline{V}_{F_N}^{(j)}. 
			$$  
		Then by adapting the proof of Theorem 7 in \citeasnoun{imbens-newey} to the triangular array for the repeated cross-section setting, for the $j$'s that satisfy Condition (iv) we obtain 
		$$
		\sqrt N {\overline{V}_{F_N}}^{- \frac 1 2}B[\hat \pi - \pi_{F_N}] \overset{F_N}{\leadsto}
		N(0, {\bf I}_{m}).
		$$
		The rest is the same as the proof of Theorem \ref{thm1}.
	\end{proof}

\section*{Appendix B: Algorithms for Computing $A$}
\label{appendix B}

This appendix details algorithms for computation of $A$. The first algorithm is the depth-first search that we in fact implemented. The second algorithm is a further refinement using Theorem \ref{prop1}. Algorithms use notation introduced in the proof of Theorem \ref{prop1}.

\medskip

\textbf{Computing A as in Theorem \ref{prop:cone}.}

\texttt{1. Initialize }$m_{1}=...=m_{J}=1$\texttt{.}

\texttt{2. Initialize }$l=2$\texttt{.}

\texttt{3. Set }$c(\mathcal{B}_{1})=x_{m_1|1},\dots,c(\mathcal{B}_{l})=x_{m_l|l}$\texttt{. Check for revealed preference cycles.}

\texttt{4. If a cycle is detected, move to step 7. Else:}

\texttt{5. If }$l<J$\texttt{, set }$l=l+1$\texttt{, }$m_{l}=1$\texttt{, and
return to step 3. Else:}

\texttt{6. Extend }$A$\texttt{\ by the column }$[m_{1},...,m_{J}]^{\prime }$%
\texttt{.}

\texttt{7a. If }$m_{l}<I_{l}$\texttt{, set }$m_{l}=m_{l}+1$\texttt{\ and
return to step 3.}

\texttt{7b. If }$m_{l}=I_{l}$\texttt{\ and }$m_{l-1}<I_{l-1}$\texttt{, set }$%
m_{l}=1$\texttt{, }$m_{l-1}=m_{l-1}+1$\texttt{, }$l=l-1$\texttt{, and return
to step 3.}

\texttt{7c. If }$m_{l}=I_{l}$\texttt{,\ }$m_{l-1}=I_{l-1}$\texttt{, and }$%
m_{l-2}<I_{l-2}$\texttt{, set }$m_{l}=m_{l-1}=1$\texttt{, }$%
m_{l-2}=m_{l-2}+1 $\texttt{, }$l=l-2$\texttt{, and return to step 3.}

\texttt{(...)}

\texttt{7z. Terminate.}

\medskip

\textbf{Refinement using Theorem \ref{prop1}}

Let budgets be arranged s.t. $(\mathcal{B}_{1},...,\mathcal{B}_{M})$ do not intersect $\mathcal{B}_{J}$;
for exposition of the algorithm, assume $\mathcal{B}_{J}$ is above these budgets.

\medskip

\texttt{1. Use preceding algorithm to compute a matrix }$A_{M+1\rightarrow J-1}$\texttt{ corresponding to budgets }$(\mathcal{B}_{M+1},...,\mathcal{B}_{J})$\texttt{, though using the full }$X$\texttt{ corresponding to budgets }$(\mathcal{B}_{1},...,\mathcal{B}_{J})$\texttt{.\footnote{This matrix has more rows than an $A$\ matrix that is only intended to apply to choice problems $(\mathcal{B}_{M+1},...,\mathcal{B}_{J})$.}}

\texttt{2. For each column }$a_{M+1\rightarrow J-1}$\texttt{\ of }$A_{M+1\rightarrow J-1}$\texttt{, go through the following steps:}

\texttt{2.1 Compute (using preceding algorithm) all vectors }$a_{1\rightarrow M}$\texttt{ s.t.}$(a_{1\rightarrow M},a_{M+1\rightarrow J-1})$\texttt{\ is rationalizable.}

\texttt{2.2 Compute (using preceding algorithm) all vectors }$a_{J}$\texttt{ s.t.}$(a_{M+1\rightarrow J-1},a_{J})$\newline \texttt{\ is rationalizable.}

\texttt{2.3 All stacked vectors }$(a_{1\rightarrow M}^{\prime},a_{M+1\rightarrow J-1}^{\prime },a_{J}^{\prime })^{\prime }$\texttt{\ are valid columns of }$A$\texttt{.}

\section*{Appendix C: Justification of Table \ref{table:sequences}}
\label{appendix C}
This appendix derives the upper bound on nodes visited by a tree search as described in Section \ref{sec:collect}. We only count nodes corresponding to $j \geq 2$ as rationalizability of implied choice behavior is checked only at those. 

Consider the number of nodes visited in generation $j+1$, i.e. corresponding to budget $\mathcal{B}_j$. Since $I_j \leq 2^{J-1}$, this is at most $2^{J-1}$ times the number of nodes in the $j$'th generation at which no choice cycle was detected. These nodes, in turn, correspond to the at most $\bar{H}_{j-1}$ direct revealed preference orderings that can occur on $(j-1)$ budgets. However, since we look at patches corresponding to the entire set of $J$ budgets, each of those orderings has multiple representations. Specifically, each patch in an $A$-matrix corresponding to the first $(j-1)$ budgets corresponds to at most $2^{J-(j-1)}$ patches in the problem under consideration (because the patches are generated by intersecting the original patch with $(J-(j-1))$ budgets). These refined patches can be arbitrarily combined across the first $(j-1)$ budgets, so that each direct revealed preference ordering on the first $(j-1)$ budgets has at most $2^{(J-j+1)(j-1)}$ representations. Thus, the number of nodes visited in generation $j+1$, $j,...,J$, is at most $\bar{H}_{j-1} 2^{(J-j+1)(j-1)+J-1}=\bar{H}_{j-1} 2^{j(J+2-j)-2}$. This bound must be summed over $j=2,...,J$.

\end{document}